

\input amstex
\expandafter\ifx\csname mathdefs.tex\endcsname\relax
  \expandafter\gdef\csname mathdefs.tex\endcsname{}
\else \message{Hey!  Apparently you were trying to
  \string twice.   This does not make sense.} 
\errmessage{Please edit your file (probably \jobname.tex) and remove
any duplicate ``\string\input'' lines} \fi




\catcode`\X=12\catcode`\@=11

\def\n@wcount{\alloc@0\count\countdef\insc@unt}
\def\n@wwrite{\alloc@7\write\chardef\sixt@@n}
\def\n@wread{\alloc@6\read\chardef\sixt@@n}
\def\r@s@t{\relax}\def\v@idline{\par}\def\@mputate#1/{#1}
\def\l@c@l#1X{\firstpart.#1}\def\gl@b@l#1X{#1}\def\t@d@l#1X{{}}

\def\crossrefs#1{\ifx\all#1\let\tr@ce=\all\else\def\tr@ce{#1,}\fi
   \n@wwrite\cit@tionsout\openout\cit@tionsout=\jobname.cit 
   \write\cit@tionsout{\tr@ce}\expandafter\setfl@gs\tr@ce,}
\def\setfl@gs#1,{\def\@{#1}\ifx\@\empty\let\next=\relax
   \else\let\next=\setfl@gs\expandafter\xdef
   \csname#1tr@cetrue\endcsname{}\fi\next}
\def\m@ketag#1#2{\expandafter\n@wcount\csname#2tagno\endcsname
     \csname#2tagno\endcsname=0\let\tail=\all\xdef\all{\tail#2,}
   \ifx#1\l@c@l\let\tail=\r@s@t\xdef\r@s@t{\csname#2tagno\endcsname=0\tail}\fi
   \expandafter\gdef\csname#2cite\endcsname##1{\expandafter
     \ifx\csname#2tag##1\endcsname\relax?\else\csname#2tag##1\endcsname\fi
     \expandafter\ifx\csname#2tr@cetrue\endcsname\relax\else
     \write\cit@tionsout{#2tag ##1 cited on page \folio.}\fi}
   \expandafter\gdef\csname#2page\endcsname##1{\expandafter
     \ifx\csname#2page##1\endcsname\relax?\else\csname#2page##1\endcsname\fi
     \expandafter\ifx\csname#2tr@cetrue\endcsname\relax\else
     \write\cit@tionsout{#2tag ##1 cited on page \folio.}\fi}
   \expandafter\gdef\csname#2tag\endcsname##1{\expandafter
      \ifx\csname#2check##1\endcsname\relax
      \expandafter\xdef\csname#2check##1\endcsname{}%
      \else\immediate\write16{Warning: #2tag ##1 used more than once.}\fi
      \multit@g{#1}{#2}##1/X%
      \write\t@gsout{#2tag ##1 assigned number \csname#2tag##1\endcsname\space
      on page \number\count0.}%
   \csname#2tag##1\endcsname}}
\def\multit@g#1#2#3/#4X{\def\t@mp{#4}\ifx\t@mp\empty%
      \global\advance\csname#2tagno\endcsname by 1 
      \expandafter\xdef\csname#2tag#3\endcsname
      {#1\number\csname#2tagno\endcsnameX}%
   \else\expandafter\ifx\csname#2last#3\endcsname\relax
      \expandafter\n@wcount\csname#2last#3\endcsname
      \global\advance\csname#2tagno\endcsname by 1 
      \expandafter\xdef\csname#2tag#3\endcsname
      {#1\number\csname#2tagno\endcsnameX}
      \write\t@gsout{#2tag #3 assigned number \csname#2tag#3\endcsname\space
      on page \number\count0.}\fi
   \global\advance\csname#2last#3\endcsname by 1
   \def\t@mp{\expandafter\xdef\csname#2tag#3/}%
   \expandafter\t@mp\@mputate#4\endcsname
   {\csname#2tag#3\endcsname\lastpart{\csname#2last#3\endcsname}}\fi}
\def\t@gs#1{\def\all{}\m@ketag#1e\m@ketag#1s\m@ketag\t@d@l p
   \m@ketag\gl@b@l r \n@wread\t@gsin
   \openin\t@gsin=\jobname.tgs \re@der \closein\t@gsin
   \n@wwrite\t@gsout\openout\t@gsout=\jobname.tgs }
\outer\def\localtags{\t@gs\l@c@l}
\outer\def\globaltags{\t@gs\gl@b@l}
\outer\def\newlocaltag#1{\m@ketag\l@c@l{#1}}
\outer\def\newglobaltag#1{\m@ketag\gl@b@l{#1}}

\newif\ifpr@ 
\def\m@kecs #1tag #2 assigned number #3 on page #4.%
   {\expandafter\gdef\csname#1tag#2\endcsname{#3}
   \expandafter\gdef\csname#1page#2\endcsname{#4}
   \ifpr@\expandafter\xdef\csname#1check#2\endcsname{}\fi}
\def\re@der{\ifeof\t@gsin\let\next=\relax\else
   \read\t@gsin to\t@gline\ifx\t@gline\v@idline\else
   \expandafter\m@kecs \t@gline\fi\let \next=\re@der\fi\next}
\def\pretags#1{\pr@true\pret@gs#1,,}
\def\pret@gs#1,{\def\@{#1}\ifx\@\empty\let\n@xtfile=\relax
   \else\let\n@xtfile=\pret@gs \openin\t@gsin=#1.tgs \message{#1} \re@der 
   \closein\t@gsin\fi \n@xtfile}

\newcount\sectno\sectno=0\newcount\subsectno\subsectno=0
\newif\ifultr@local \def\ultralocal{\ultr@localtrue}
\def\firstpart{\number\sectno}
\def\lastpart#1{\ifcase#1 \or a\or b\or c\or d\or e\or f\or g\or h\or 
   i\or k\or l\or m\or n\or o\or p\or q\or r\or s\or t\or u\or v\or w\or 
   x\or y\or z \fi}

\def\resetall{\global\advance\sectno by 1\subsectno=0
   \gdef\firstpart{\number\sectno}\r@s@t}
\def\resetsub{\global\advance\subsectno by 1
   \gdef\firstpart{\number\sectno.\number\subsectno}\r@s@t}
\def\newsection#1\par{\resetall\vskip0pt plus.3\vsize\penalty-250
   \vskip0pt plus-.3\vsize\bigskip\bigskip
   \message{#1}\leftline{\bf#1}\nobreak\bigskip}
\def\subsection#1\par{\ifultr@local\resetsub\fi
   \vskip0pt plus.2\vsize\penalty-250\vskip0pt plus-.2\vsize
   \bigskip\smallskip\message{#1}\leftline{\bf#1}\nobreak\medskip}

\def\t@gsoff#1,{\def\@{#1}\ifx\@\empty\let\next=\relax\else\let\next=\t@gsoff
   \def\@@{p}\ifx\@\@@\else
   \expandafter\gdef\csname#1cite\endcsname##1{\zeigen{##1}}
   \expandafter\gdef\csname#1page\endcsname##1{?}
   \expandafter\gdef\csname#1tag\endcsname##1{\zeigen{##1}}\fi\fi\next}
\def\verbatimtags{\ifx\all\relax\else\expandafter\t@gsoff\all,\fi}
\def\zeigen#1{\hbox{$\langle$}#1\hbox{$\rangle$}}

\def\(#1){\edef\dot@g{\ifmmode\ifinner(\hbox{\noexpand\etag{#1}})
   \else\noexpand\eqno(\hbox{\noexpand\etag{#1}})\fi
   \else(\noexpand\ecite{#1})\fi}\dot@g}

\newif\ifbr@ck
\def\eat#1{}
\def\[#1]{\br@cktrue[\br@cket#1'X]}
\def\br@cket#1'#2X{\def\temp{#2}\ifx\temp\empty\let\next\eat
   \else\let\next\br@cket\fi
   \ifbr@ck\br@ckfalse\br@ck@t#1,X\else\br@cktrue#1\fi\next#2X}
\def\br@ck@t#1,#2X{\def\temp{#2}\ifx\temp\empty\let\neext\eat
   \else\let\neext\br@ck@t\def\temp{,}\fi
   \def\teemp{#1}\ifx\teemp\empty\else\rcite{#1}\fi\temp\neext#2X}
\def\resetbr@cket{\gdef\[##1]{[\rtag{##1}]}}
\def\references{\resetbr@cket\newsection References\par}

\newtoks\symb@ls\newtoks\s@mb@ls\newtoks\p@gelist\n@wcount\ftn@mber
    \ftn@mber=1\newif\ifftn@mbers\ftn@mbersfalse\newif\ifbyp@ge\byp@gefalse
\def\defm@rk{\ifftn@mbers\n@mberm@rk\else\symb@lm@rk\fi}
\def\n@mberm@rk{\xdef\m@rk{{\the\ftn@mber}}%
    \global\advance\ftn@mber by 1 }
\def\rot@te#1{\let\temp=#1\global#1=\expandafter\r@t@te\the\temp,X}
\def\r@t@te#1,#2X{{#2#1}\xdef\m@rk{{#1}}}
\def\b@@st#1{{$^{#1}$}}\def\str@p#1{#1}
\def\symb@lm@rk{\ifbyp@ge\rot@te\p@gelist\ifnum\expandafter\str@p\m@rk=1 
    \s@mb@ls=\symb@ls\fi\write\f@nsout{\number\count0}\fi \rot@te\s@mb@ls}
\def\byp@ge{\byp@getrue\n@wwrite\f@nsin\openin\f@nsin=\jobname.fns 
    \n@wcount\currentp@ge\currentp@ge=0\p@gelist={0}
    \re@dfns\closein\f@nsin\rot@te\p@gelist
    \n@wread\f@nsout\openout\f@nsout=\jobname.fns }
\def\m@kelist#1X#2{{#1,#2}}
\def\re@dfns{\ifeof\f@nsin\let\next=\relax\else\read\f@nsin to \f@nline
    \ifx\f@nline\v@idline\else\let\t@mplist=\p@gelist
    \ifnum\currentp@ge=\f@nline
    \global\p@gelist=\expandafter\m@kelist\the\t@mplistX0
    \else\currentp@ge=\f@nline
    \global\p@gelist=\expandafter\m@kelist\the\t@mplistX1\fi\fi
    \let\next=\re@dfns\fi\next}
\def\symbols#1{\symb@ls={#1}\s@mb@ls=\symb@ls} 
\def\bigsymbol{\textstyle}
\symbols{\bigsymbol\ast,\dagger,\ddagger,\sharp,\flat,\natural,\star}
\def\ftnumbers{\ftn@mberstrue} \def\ftsymbols{\ftn@mbersfalse}
\def\paginal{\byp@ge} \def\resetftnumbers{\ftn@mber=1}
\def\ftnote#1{\defm@rk\expandafter\expandafter\expandafter\footnote
    \expandafter\b@@st\m@rk{#1}}

\long\def\jump#1\endjump{}
\def\ssum{\mathop{\lower .1em\hbox{$\textstyle\Sigma$}}\nolimits}

\def\qed{\nobreak\kern 1em \vrule height .5em width .5em depth 0em}
\def\newneq{\hbox{\rlap{\hbox to 1\wd9{\hss$=$\hss}}\raise .1em 
   \hbox to 1\wd9{\hss$\scriptscriptstyle/$\hss}}}
\def\subsetne{\setbox9 = \hbox{$\subset$}\mathrel{\hbox{\rlap
   {\lower .4em \newneq}\raise .13em \hbox{$\subset$}}}}
\def\supsetne{\setbox9 = \hbox{$\subset$}\mathrel{\hbox{\rlap
   {\lower .4em \newneq}\raise .13em \hbox{$\supset$}}}}

\def\vbar{\mathchoice{\vrule height6.3ptdepth-.5ptwidth.8pt\kern-.8pt}
   {\vrule height6.3ptdepth-.5ptwidth.8pt\kern-.8pt}
   {\vrule height4.1ptdepth-.35ptwidth.6pt\kern-.6pt}
   {\vrule height3.1ptdepth-.25ptwidth.5pt\kern-.5pt}}
\def\f@dge{\mathchoice{}{}{\mkern.5mu}{\mkern.8mu}}
\def\b@c#1#2{{\rm \mkern#2mu\vbar\mkern-#2mu#1}}
\def\b@b#1{{\rm I\mkern-3.5mu #1}}
\def\b@a#1#2{{\rm #1\mkern-#2mu\f@dge #1}}
\def\bb#1{{\count4=`#1 \advance\count4by-64 \ifcase\count4\or\b@a A{11.5}\or
   \b@b B\or\b@c C{5}\or\b@b D\or\b@b E\or\b@b F \or\b@c G{5}\or\b@b H\or
   \b@b I\or\b@c J{3}\or\b@b K\or\b@b L \or\b@b M\or\b@b N\or\b@c O{5} \or
   \b@b P\or\b@c Q{5}\or\b@b R\or\b@a S{8}\or\b@a T{10.5}\or\b@c U{5}\or
   \b@a V{12}\or\b@a W{16.5}\or\b@a X{11}\or\b@a Y{11.7}\or\b@a Z{7.5}\fi}}

\catcode`\X=11 \catcode`\@=12

\expandafter\ifx\csname citeadd.tex\endcsname\relax
\expandafter\gdef\csname citeadd.tex\endcsname{}
\else \message{Hey!  Apparently you were trying to
\string twice.   This does not make sense.} 
\errmessage{Please edit your file (probably \jobname.tex) and remove
any duplicate ``\string\input'' lines} \fi

\sectno=-2   
\localtags
\NoBlackBoxes
\ifx\shlhetal\undefinedcontrolsequence\let\shlhetal\relax\fi
\define\mr{\medskip\roster}
\define\sn{\smallskip\noindent}
\define\mn{\medskip\noindent}
\define\bn{\bigskip\noindent}
\define\ub{\underbar}
\define\wilog{\text{without loss of generality}}
\define\ermn{\endroster\medskip\noindent}

\define\dbcu{\dsize\bigcup}
\define \nl{\newline}
\documentstyle {amsppt}
\topmatter
\title{{\it THE GENERALIZED CONTINUUM HYPOTHESIS REVISITED}} \endtitle
\author {Saharon Shelah \thanks {\null\newline Partially supported by the 
Israeli Basic Research Fund and the BSF and I would like to thank Alice 
Leonhardt for the beautiful typing.\newline
 Done - 8/91(\S1), 9/91(\S2,\S3), \newline
 Latest Revision - 98/Sept/19 \newline
 Publ. No.460} \endthanks} \endauthor
\affil{Institute of Mathematics\\
 The Hebrew University\\
 Jerusalem, Israel
 \medskip
 Rutgers University\\
 Mathematics Department\\
 New Brunswick, NJ  USA} \endaffil
\endtopmatter
\document  

\expandafter\ifx\csname alice2jlem.tex\endcsname\relax
  \expandafter\gdef\csname alice2jlem.tex\endcsname{}
\else \message{Hey!  Apparently you were trying to
\string  twice.   This does not make sense.}
\errmessage{Please edit your file (probably \jobname.tex) and remove
any duplicate ``\string\input'' lines} \fi

\expandafter\ifx\csname bib4plain.tex\endcsname\relax
  \expandafter\gdef\csname bib4plain.tex\endcsname{}
\else \message{Hey!  Apparently you were trying to \string twice.   This does not make sense.}
\errmessage{Please edit your file (probably \jobname.tex) and remove
any duplicate ``\string\input'' lines} \fi

\def\renewcommand{\newcommand}	       
\edef\cite{\the\catcode`@}%
\catcode`@ = 11
\let\@oldatcatcode = \cite
\chardef\@letter = 11
\chardef\@other = 12
%
%
%
%
\def\@innerdef#1#2{\edef#1{\expandafter\noexpand\csname #2\endcsname}}%
%
%
\@innerdef\@innernewcount{newcount}%
\@innerdef\@innernewdimen{newdimen}%
\@innerdef\@innernewif{newif}%
\@innerdef\@innernewwrite{newwrite}%
%
%
%
\def\@gobble#1{}%
%
%
%
\ifx\inputlineno\@undefined
   \let\@linenumber = \empty 
\else
   \def\@linenumber{\the\inputlineno:\space}%
\fi
%
%
%
\def\@futurenonspacelet#1{\def\cs{#1}%
   \afterassignment\@stepone\let\@nexttoken=
}%
\begingroup 
\def\\{\global\let\@stoken= }%
\\ 
\endgroup
\def\@stepone{\expandafter\futurelet\cs\@steptwo}%
\def\@steptwo{\expandafter\ifx\cs\@stoken\let\@@next=\@stepthree
   \else\let\@@next=\@nexttoken\fi \@@next}%
\def\@stepthree{\afterassignment\@stepone\let\@@next= }%
%
%
%
\def\@getoptionalarg#1{%
   \let\@optionaltemp = #1%
   \let\@optionalnext = \relax
   \@futurenonspacelet\@optionalnext\@bracketcheck
}%
%
%
\def\@bracketcheck{%
   \ifx [\@optionalnext
      \expandafter\@@getoptionalarg
   \else
      \let\@optionalarg = \empty
      \expandafter\@optionaltemp
   \fi
}%
\def\@@getoptionalarg[#1]{%
   \def\@optionalarg{#1}%
   \@optionaltemp
}%
%
%
%
\def\@nnil{\@nil}%
\def\@fornoop#1\@@#2#3{}%
\def\@for#1:=#2\do#3{%
   \edef\@fortmp{#2}%
   \ifx\@fortmp\empty \else
      \expandafter\@forloop#2,\@nil,\@nil\@@#1{#3}%
   \fi
}%
\def\@forloop#1,#2,#3\@@#4#5{\def#4{#1}\ifx #4\@nnil \else
       #5\def#4{#2}\ifx #4\@nnil \else#5\@iforloop #3\@@#4{#5}\fi\fi
}%
\def\@iforloop#1,#2\@@#3#4{\def#3{#1}\ifx #3\@nnil
       \let\@nextwhile=\@fornoop \else
      #4\relax\let\@nextwhile=\@iforloop\fi\@nextwhile#2\@@#3{#4}%
}%
%
%
%
\@innernewif\if@fileexists
\def\@testfileexistence{\@getoptionalarg\@finishtestfileexistence}%
\def\@finishtestfileexistence#1{%
   \begingroup
      \def\extension{#1}%
      \immediate\openin0 =
         \ifx\@optionalarg\empty\jobname\else\@optionalarg\fi
         \ifx\extension\empty \else .#1\fi
         \space
      \ifeof 0
         \global\@fileexistsfalse
      \else
         \global\@fileexiststrue
      \fi
      \immediate\closein0
   \endgroup
}%
%
%
%
%
\def\bibliographystyle#1{%
   \@readauxfile
   \@writeaux{\string\bibstyle{#1}}%
}%
\let\bibstyle = \@gobble
%
%
\let\bblfilebasename = \jobname
\def\bibliography#1{%
   \@readauxfile
   \@writeaux{\string\bibdata{#1}}%
   \@testfileexistence[\bblfilebasename]{bbl}%
   \if@fileexists
      \nobreak
      \@readbblfile
   \fi
}%
\let\bibdata = \@gobble
%
%
\def\nocite#1{%
   \@readauxfile
   \@writeaux{\string\citation{#1}}%
}%
\@innernewif\if@notfirstcitation
%
%
\def\cite{\@getoptionalarg\@cite}%
%
%
\def\@cite#1{%
   \let\@citenotetext = \@optionalarg
   \printcitestart
   \nocite{#1}%
   \@notfirstcitationfalse
   \@for \@citation :=#1\do
   {%
      \expandafter\@onecitation\@citation\@@
   }%
   \ifx\empty\@citenotetext\else
      \printcitenote{\@citenotetext}%
   \fi
   \printcitefinish
}%
\def\@onecitation#1\@@{%
   \if@notfirstcitation
      \printbetweencitations
   \fi
   \expandafter \ifx \csname\@citelabel{#1}\endcsname \relax
      \if@citewarning
         \message{\@linenumber Undefined citation `#1'.}%
      \fi
      \expandafter\gdef\csname\@citelabel{#1}\endcsname{%
\strut
\vadjust{\vskip-\dp\strutbox
\vbox to 0pt{\vss\parindent0cm \leftskip=\hsize 
\advance\leftskip3mm
\advance\hsize 4cm\strut\openup-4pt 
\rightskip 0cm plus 1cm minus 0.5cm ?  #1 ?\strut}}
         {\tt
            \escapechar = -1
            \nobreak\hskip0pt
            \expandafter\string\csname#1\endcsname
            \nobreak\hskip0pt
         }%
      }%
   \fi
   \csname\@citelabel{#1}\endcsname
   \@notfirstcitationtrue
}%
%
%
\def\@citelabel#1{b@#1}%
%
%
\def\@citedef#1#2{\expandafter\gdef\csname\@citelabel{#1}\endcsname{#2}}%
%
%
%
\def\@readbblfile{%
   \ifx\@itemnum\@undefined
      \@innernewcount\@itemnum
   \fi
   \begingroup
      \def\begin##1##2{%
         \setbox0 = \hbox{\biblabelcontents{##2}}%
         \biblabelwidth = \wd0
      }%
      \def\end##1{}
      %
      %
      \@itemnum = 0
      \def\bibitem{\@getoptionalarg\@bibitem}%
      \def\@bibitem{%
         \ifx\@optionalarg\empty
            \expandafter\@numberedbibitem
         \else
            \expandafter\@alphabibitem
         \fi
      }%
      \def\@alphabibitem##1{%
         \expandafter \xdef\csname\@citelabel{##1}\endcsname {\@optionalarg}%
         \ifx\biblabelprecontents\@undefined
            \let\biblabelprecontents = \relax
         \fi
         \ifx\biblabelpostcontents\@undefined
            \let\biblabelpostcontents = \hss
         \fi
         \@finishbibitem{##1}%
      }%
      \def\@numberedbibitem##1{%
         \advance\@itemnum by 1
         \expandafter \xdef\csname\@citelabel{##1}\endcsname{\number\@itemnum}%
         \ifx\biblabelprecontents\@undefined
            \let\biblabelprecontents = \hss
         \fi
         \ifx\biblabelpostcontents\@undefined
            \let\biblabelpostcontents = \relax
         \fi
         \@finishbibitem{##1}%
      }%
      \def\@finishbibitem##1{%
         \biblabelprint{\csname\@citelabel{##1}\endcsname}%
         \@writeaux{\string\@citedef{##1}{\csname\@citelabel{##1}\endcsname}}%
         \ignorespaces
      }%
      %
      %
      \let\em = \bblem
      \let\newblock = \bblnewblock
      \let\sc = \bblsc
      \frenchspacing
      \clubpenalty = 4000 \widowpenalty = 4000
      \tolerance = 10000 \hfuzz = .5pt
      \everypar = {\hangindent = \biblabelwidth
                      \advance\hangindent by \biblabelextraspace}%
      \bblrm
      \parskip = 1.5ex plus .5ex minus .5ex
      \biblabelextraspace = .5em
      \bblhook
      \input \bblfilebasename.bbl
   \endgroup
}%
%
%
\@innernewdimen\biblabelwidth
\@innernewdimen\biblabelextraspace
%
%
%
\def\biblabelprint#1{%
   \noindent
   \hbox to \biblabelwidth{%
      \biblabelprecontents
      \biblabelcontents{#1}%
      \biblabelpostcontents
   }%
   \kern\biblabelextraspace
}%
%
%
%
\def\biblabelcontents#1{{\bblrm [#1]}}%
%
%
\def\bblrm{\rm}%
%
%
\def\bblem{\it}%
%
%
\def\bblsc{\ifx\@scfont\@undefined
              \font\@scfont = cmcsc10
           \fi
           \@scfont
}%
%
%
\def\bblnewblock{\hskip .11em plus .33em minus .07em }%
%
%
\let\bblhook = \empty
%
%
%
\def\printcitestart{[}
\def\printcitefinish{]}
\def\printbetweencitations{, }
\def\printcitenote#1{, #1}
%
%
%
\let\citation = \@gobble
%
%
%
\@innernewcount\@numparams
%
%
\def\newcommand#1{%
   \def\@commandname{#1}%
   \@getoptionalarg\@continuenewcommand
}%
%
%
\def\@continuenewcommand{%
   \@numparams = \ifx\@optionalarg\empty 0\else\@optionalarg \fi \relax
   \@newcommand
}%
%
%
\def\@newcommand#1{%
   \def\@startdef{\expandafter\edef\@commandname}%
   \ifnum\@numparams=0
      \let\@paramdef = \empty
   \else
      \ifnum\@numparams>9
         \errmessage{\the\@numparams\space is too many parameters}%
      \else
         \ifnum\@numparams<0
            \errmessage{\the\@numparams\space is too few parameters}%
         \else
            \edef\@paramdef{%
               \ifcase\@numparams
                  \empty  No arguments.
               \or ####1%
               \or ####1####2%
               \or ####1####2####3%
               \or ####1####2####3####4%
               \or ####1####2####3####4####5%
               \or ####1####2####3####4####5####6%
               \or ####1####2####3####4####5####6####7%
               \or ####1####2####3####4####5####6####7####8%
               \or ####1####2####3####4####5####6####7####8####9%
               \fi
            }%
         \fi
      \fi
   \fi
   \expandafter\@startdef\@paramdef{#1}%
}%
%
%
%
%
\def\@readauxfile{%
   \if@auxfiledone \else 
      \global\@auxfiledonetrue
      \@testfileexistence{aux}%
      \if@fileexists
         \begingroup
            \endlinechar = -1
            \catcode`@ = 11
            \input \jobname.aux
         \endgroup
      \else
         \message{\@undefinedmessage}%
         \global\@citewarningfalse
      \fi
      \immediate\openout\@auxfile = \jobname.aux
   \fi
}%
%
%
\newif\if@auxfiledone
\ifx\noauxfile\@undefined \else \@auxfiledonetrue\fi
%
%
%
%
\@innernewwrite\@auxfile
\def\@writeaux#1{\ifx\noauxfile\@undefined \write\@auxfile{#1}\fi}%
%
%
%
\ifx\@undefinedmessage\@undefined
   \def\@undefinedmessage{No .aux file; I won't give you warnings about
                          undefined citations.}%
\fi
%
%
\@innernewif\if@citewarning
\ifx\noauxfile\@undefined \@citewarningtrue\fi
%
%
%
\catcode`@ = \@oldatcatcode


\def\widestnumber#1#2{}

\def\rm{\fam0 \tenrm}

\def\fakesubhead#1\endsubhead{\bigskip\noindent{\bf#1}\par}


%
%
%

%

\font\textrsfs=rsfs10
\font\scriptrsfs=rsfs7
\font\scriptscriptrsfs=rsfs5

\newfam\rsfsfam
\textfont\rsfsfam=\textrsfs
\scriptfont\rsfsfam=\scriptrsfs
\scriptscriptfont\rsfsfam=\scriptscriptrsfs

\edef\oldcatcodeofat{\the\catcode`\@}
\catcode`\@11

\def\Cal@@#1{\noaccents@ \fam \rsfsfam #1}

\catcode`\@\oldcatcodeofat

\newpage

\head {Annotated Content} \endhead  \resetall \bn
\S0 Introduction
\mr
\item "{{}}" [We explain why we consider the main theorem here a
reasonable revision of GCH (but provable in ZFC).]
\ermn
\mn
\S1 The generic ultrapower proof
\mr
\item "{{}}"  [We prove that for $\mu$ strong limit $> \aleph_0$ for every
$\lambda > \mu$ for some $\kappa < \mu$, there is ${\Cal P} \subseteq
[\lambda]^{< \mu}$ of cardinality $\lambda$ such that every $A \in
[\lambda]^{< \mu}$ is the union of $< \kappa$ members of $P$.  We do it using
generic ultrapowers.  We draw some immediate conclusions.]
\ermn
\mn
\S2 The main theorem revisited
\mr
\item "{{}}"  [We give a somewhat stronger theorem, with a proof using pcf
calculus.   We then show that if $\mu$ is a limit cardinal satisfying the
conclusion of the main theorem, then ${\frak a} \subseteq \text{ Reg}
\backslash \mu \and |{\frak a}| < \mu \Rightarrow |\text{pcf}({\frak a})|
\le \mu$.  Then we prove a converse: if $(\forall {\frak a} \subseteq
\text{ Reg} \backslash \mu)(|{\frak a}| < \mu \Rightarrow |\text{pcf}
({\frak a})| < \mu)$ (or somewhat less) then the conclusion on the main 
theorem holds.]
\ermn
\mn
\S3  Application
\mr
\item "{{}}"  [We draw a conclusion on diamonds (and $(D \ell)$), hence on
the omitting types theorem for $L(Q)$ in the $\lambda^+$-interpretation and
on relatives of Arhangelskii's problem.]
\ermn
\mn
\ub{Appendix: Existence of tiny models}
\mr
\item "{{}}"  [We show the close connection of the problems to the existence 
of families of large sets with pairwise finite intersections.]
\endroster
\newpage

\head {\S0 Introduction} \endhead  \resetall \bigskip

I had a dream, quite a natural one for a mathematician in the
twentieth century: to solve a Hilbert problem, preferably positively.
This is quite hard for (at least) three reasons:
\roster
\item "{(a)}"  those problems are almost always hard
\sn
\item "{(b)}"  almost all have been solved
\sn
\item "{(c)}"  my (lack of) knowledge excludes almost all.
\ermn
Now (c) points out the first Hilbert problem as it is in set
theory; also being the first it occupy a place of honor.\newline

The problem asks ``is the continuum hypothesis true?'', i.e.,
\medskip
\roster
\item  is $2^{\aleph_0} = \aleph_1$?
\endroster
\medskip
\noindent
More generally, is the generalized continuum hypothesis true?  Which means:
\medskip
\roster
\item "{(2)}"  is $2^{\aleph_\alpha} = \aleph_{\alpha+1}$ for all
ordinals $\alpha$?
\endroster
\medskip
\noindent
I think the meaning of the question is what are the laws of cardinal
arithmetic; it was known that addition and multiplication of infinite
cardinals
is ``trivial'', i.e. previous generations have not left us anything to solve:
$$
\lambda + \mu = \lambda \times \mu = \text{ max}\{ \lambda,\mu \}.
$$
This would have certainly made elementary school pupils happier than the usual
laws, but we have been left with exponentiation only.  As there were two
operations on infinite cardinals increasing them --- $2^\lambda$ and
$\lambda^+$ --- it was most natural to assume that those two operations 
are the same; in fact, in this case also exponentiation becomes very 
simple; usually $\lambda^\mu = \text{ max} \{ \lambda,\mu^+ \}$, the exception
being that when cf$(\lambda) \le \mu < \lambda$ we have 
$\lambda^\mu = \lambda^+$ where \newline
cf$(\lambda) =: \text{ min}\{ \kappa:\text{ there are } \lambda_i < 
\lambda \text{ for } i < \kappa \text{ such that } \lambda = 
\dsize \sum_{i<\kappa}\lambda_i \}$.
Non-set theorists may be reminded that $\lambda = \mu^+$ if $\mu =
\aleph_\alpha$ and $\lambda = \aleph_{\alpha +1}$, and then 
$\lambda$ is called the successor of $\mu$ and we know 
cf$(\aleph_{\alpha +1}) = \aleph
_{\alpha +1}$; we call a cardinal $\lambda$ regular if cf$(\lambda) = \lambda$
and singular otherwise.  So successor cardinals are regular and also
$\aleph_0$, but it is ``hard to come by" other regular cardinals so we may
ignore them.  Note $\aleph_\omega = \dsize \sum_{n < \omega} \aleph_n$ is the
first singular cardinal, and for $\delta$ a limit ordinal $> |\delta|$ we have
$\aleph_\delta$ singular, but there are limit $\delta = \aleph_\delta$ for
which $\aleph_\delta$ is singular.

Probably the interpretation of Hilbert's first problem as ``find all laws of
cardinal arithmetic'' is too broad \footnote{On this see \cite{Sh:g} or
\cite{Sh:400a}, note that under the interpretation of the problem 
there is much to say.}, still ``is cardinal arithmetic simple'' is
a reasonable interpretation.

Unfortunately, there are some ``difficulties".
On the one hand, G\"odel 
had proved that GCH may be true (specifically it holds
in the universe of constructible sets, called $L$).  
On the other hand, Cohen had proved that CH may be false 
(by increasing the universe of sets by forcing), in
fact, $2^{\aleph_0}$ can be anything reasonable, i.e., cf$(2^{\aleph_0}) > 
\aleph_0$. 

Continuing Cohen, Solovay proved that $2^{\aleph_n}$ for $n < \omega$ can be
anything reasonable: it should be non-decreasing and cf$(2^\lambda) >
\lambda$.  Continuing this, Easton proved that the function $\lambda
\mapsto 2^\lambda$ for regular cardinals is arbitrary (except for the laws
above).  Well, we can still hope to salvage something by proving that 
(2) holds for
``most'' cardinals; unfortunately, Magidor had proved the consistency of
$2^\lambda > \lambda^+$ for all $\lambda$ in any pregiven initial segment of
the cardinals and then Foreman and Woodin \cite{FW} for all $\lambda$.

Such difficulties should not deter the truly dedicated ones; first note
that we should not identify exponentiation with the specific case of
exponentiation $2^\lambda$, in fact Easton's results indicate that
on this (for $\lambda$ regular) we cannot say anything more, but they do not 
rule out saying something on
$\lambda^\mu$ when $\mu < \lambda$, and we can rephrase the GCH as
\roster
\item "{(3)}"  for every regular $\kappa < \lambda$ we have $\lambda^\kappa
= \lambda$.
\endroster
\medskip
Ahah, now that we have two parameters we can
look again at ``for most pairs of cardinals (3) holds.''  However, this is a
bad division, because, say, a failure for $\kappa = \aleph_1$ implies a
failure for $\kappa = \aleph_0$. \newline
To rectify this we suggest another division, we
define ``$\lambda$ to the revised power of $\kappa$", for $\kappa$ regular
$< \lambda$
as
$$
\align
\lambda^{[\kappa]} = \text{ Min} \biggl\{
|{\Cal P}|:&{\Cal P} \text{ a family of subsets of } \lambda
\text{ each of cardinality } \kappa \\
  &\text{ such that any subset of } \lambda \text{ of cardinality } \kappa \\
  &\text{ is contained in the union of } < \kappa
\text{ members of } {\Cal P} \biggr\}.
\endalign
$$

\noindent
This answers the criticism above and is a better slicing because: \newline
\medskip
\roster
\item "{(A)}" $\text{ for every } \lambda > \kappa
 \text{ we have: } \lambda^\kappa =
\lambda \text{ \underbar{iff} } 2^\kappa \le \lambda
\text{ and for every regular } \theta \le \kappa$, \newline
$\lambda^{[\theta]} = \lambda$.
\item "{(B)}"  By Gitik, Shelah \cite{GiSh:344}, the values of, e.g.,
$\lambda^{[\aleph_0]},\dotsc,\lambda^{[\aleph_n]}$ are essentially
independent.
\endroster
\medskip

\noindent
Now we rephrase the generalized continuum hypothesis as:
\medskip
\roster
\item "{$(4)$}"  $\text{ for most pairs } (\lambda,\kappa),
 \,\, \lambda^{[\kappa]} = \lambda$
\endroster
\medskip

\noindent
Is such reformulation legitimate?  As an argument, I can cite, from the
book \cite{Br} on Hilbert's problems, Lorentz's article on the thirteenth 
problem.  The problem
was
\medskip
\roster
\item "{$(*)$}"  Prove that the equation of the seventh degree $x^7 + ax^3
+ bx^2 + cx + 1 = 0$ is not solvable with the help of any continuous functions
of only two variables.
\endroster
\medskip

Lorentz does not even discuss the change from $7$ to $n$ and he shortly
changes it to (see \cite[Ch.II,p.419]{Br})
\medskip
\roster
\item "{$(*)'$}"  Prove that there are continuous functions of three variables
not represented by continuous functions of two variables.
\endroster
\medskip

\noindent
Then, he discusses Kolmogorov's solution and improvements.  He opens
the second section with (\cite[p.421,16-22]{Br}): ``that having disproved 
the conjecture is not
solving it, we should reformulate the problem in the light of the
counterexamples and prove it, which in his case: (due to Vituvskin) the
fundamental theorem of the Differential Calculus: there are $r$-times
continuously differential functions of $n$ variables not represented by
superpositions of $r$ times continuously times differential functions of less
than $n$ variables".

Concerning the fifth problem, Gleason (who makes a major contribution to its
solution) says (in \cite{AAC90}): ``Of course, many mathematicians are not 
aware that the
problem as stated by Hilbert is not the problem that has been ultimately
called the Fifth Problem.  It was shown very, very early that what he was
asking people to consider was actually false.  He asked to show that the
action of a locally-euclidean group on a manifold was always analytic, and
that's false.  It's only the group itself that's analytic, the action on a
manifold need not be.  So you had to change things considerably before you
could make the statement he was concerned with true.  That's sort of
interesting, I think.  It's also part of the way a mathematical theory
develops.  People have ideas about what ought to be so and they propose this
as a good question to work on, and then it turns out that part of it isn't
so."

In our case, I feel that while the discovery of $L$ (the constructible
universe) by G\"odel and the discovery of forcing by Cohen are fundamental
discoveries in set theory, things which
are and will continue to be in its center, forming a basis for flourishing
research, and they provide for the first Hilbert problem
a negative solution which justifies our reinterpretation of it.  Of course,
it is very reasonable to include independence results in a reinterpretation.

Back to firmer grounds, how will we interpret ``for most''? 
The simplest ways are to say ``for each
$\lambda$ for most $\kappa$'' or ``for each $\kappa$ for most $\lambda$''.
The second interpretation holds in a non-interesting way: for each $\kappa$
for many $\lambda$'s, $\lambda^\kappa = \lambda$ hence $\lambda^{[\kappa]}
= \lambda$ (e.g. $\mu^\kappa$ when $\mu \ge 2$).  So the best we can hope for
is: for every $\lambda$ for most small $\kappa$'s (remember we have restricted
ourselves to regular $\kappa$ quite smaller than $\lambda$).  To fix the
difference we restrict ourselves to $\lambda > \beth_\omega > \kappa$.  Now
what is a reasonable interpretation of ``for most $\kappa < \beth_\omega$''?
The reader may well stop and reflect.  As ``all is forbidden'' 
(by \cite{GiSh:344} even finitely many exceptions are possible), 
the simplest offer I think is ``for all but boundedly many". \newline
\noindent
So the best we can hope for is ($\beth_\omega$ is for definiteness):
\medskip
\roster
\item "{$(5)$}"  $\text{ if } \lambda > \beth_\omega,
\text{ for every large enough regular }
\kappa < \beth_\omega,\lambda^{[\kappa]} = \lambda$ \nl
(and similarly replacing $\beth_\omega$ by any strong limit cardinal).
\endroster
\medskip

If the reader has agreed so far, he is trapped into admitting that here we
solved Hilbert's first problem positively (see \scite{0.1} below).  
Now we turn from fun to business.
\mn
A consequence is
\mr
\item "{$(*)_6$}"  for every $\lambda \ge \beth_\omega$ for some $n$ and
\footnote{where $[\lambda]^{< \kappa} = \{a \subseteq \lambda:|a| <
\kappa\}$} ${\Cal P} \subseteq [\lambda]^{< \beth_\omega}$ of cardinality
$\lambda$, every $a \in [\lambda]^{< \beth_\omega}$ is the union of
$< \beth_n$ members of ${\Cal P}$.
\ermn
The history above was written just to lead to (5), 
for a fuller history see \cite{Sh:g}. 
\bn
More fully our main result is
\proclaim{\stag{0.1} The revised GCH theorem}  Assume we fix an uncountable
strong limit cardinal $\mu$ (i.e., $\mu > \aleph_0,(\forall \theta < \mu)
(2^\theta < \mu)$, e.g. $\mu = \beth_\omega = \sum \beth_n$ where $\beth_0 =
\aleph_0,\beth_{n+1} = 2^{\beth_n}$). \nl
\ub{Then} for every $\lambda \ge \mu$ for some $\kappa < \mu$ we have:
\mr
\item "{$(a)$}"  $\kappa \le \theta < \mu \Rightarrow \lambda^{[\theta]}
= \lambda$
\sn
\item "{$(b)$}"  there is a family ${\Cal P}$ of $\lambda$ subsets of
$\lambda$ each of cardinality $< \mu$ such that every subset of $\lambda$ of
cardinality $\mu$ is equal to the union of $< \kappa$ members of
${\Cal P}$.
\endroster
\endproclaim
\bigskip

\demo{Proof}  It is enough to prove it for singular $\mu$.

Clause (a) follows by clause (b) (just use ${\Cal P}_\theta
= \{a \in {\Cal P}:|a| \le \theta\}$) and clause (b) holds by
\scite{1.2}(4).
\enddemo
\bigskip

In \S1 we prove the theorem using a generic embedding based on 
\cite[Ch.VI,\S1]{Sh:g} (hence using simple forcing) and give some 
applications, mainly, they are 
reformulations.  For example, for $\lambda \ge \beth_\omega$ for every
$\theta < \beth_\omega$ large
enough, there is no tree with $\lambda$ nodes and $> \lambda$ 
$\theta$-branches.  Also we explain that this is sufficient for proving that
e.g. a topology (not necessarily even $T_0$!) with a base of cardinality 
$\mu \ge \beth_\omega$ and $> \mu$ open sets has at least 
$\beth_{\omega + 1}$ open sets.

In \scite{2.1} we give another proof (so not relying on \S1),
more inside pcf theory and saying somewhat more.  In \scite{2.6} 
we show that a
property of $\mu = \beth_\omega$ which suffices is: $\mu$ is a limit cardinal
such that $|{\frak a}| < \mu \Rightarrow |\text{pcf}({\frak a})| < \mu$ 
giving a third proof.   This is almost a converse to \scite{2.5}.
Now \S3 deals with applications: we show that for $\lambda \ge \beth_\omega$,
$2^\lambda = \lambda^+$ is equivalent to $\diamondsuit_{\lambda^+}$ (moreover
$\lambda = \lambda^{<\lambda}$ is equivalent to $(D\ell)_\lambda$, a weak
version of diamond).  We also deal with a general topology problem: can
every space be divided to two pieces, no one containing a compactum (say a
topological copy
of ${}^\omega 2$), showing its connection to pcf theory, and proving a
generalization when the cardinal parameter is $> \beth_\omega$. 
Lastly we prove there are no tiny models for theories with a non-trivial
type (see \cite{LaPiRo}) of cardinality $\ge \beth_\omega$, partially 
solving a problem from Laskowski, Pillay and Rothmaler \cite{LaPiRo}.
\sn
For other applications see \cite[\S8]{Sh:575}.  This work is continued in
\cite{Sh:513}, for further discussion see \cite{Sh:666}.  For more on
Arhangelskii's problem see \cite{Sh:668}.

We thank Todd Eisworth for many corrections and improving presentation.
\newpage 

\head {\S1 The generic ultrapower proof} \endhead  \resetall 
\proclaim{\stag{1.1} Theorem}  Assume  $\mu$  is strong limit singular and
$\lambda  > \mu$. \underbar{Then} there are only boundedly many
$\kappa  < \mu$  such that for some $\theta \in (\mu,\lambda)$ we have
{\rm pp}$_{\Gamma(\mu^+,\kappa)}(\theta) \ge \lambda$ (so $\kappa \le 
\text{{\rm cf\/}}(\theta) < \mu < \theta)$.
\endproclaim
\bn
We list some conclusions, which are immediate by older works.
\demo{\stag{1.2} Conclusion}  For every 
$\mu$ strong limit such that cf$(\mu) = 
\sigma < \mu < \lambda$, for some  $\kappa < \mu$ we have:
\roster
\item  for every  ${\frak a} \subseteq \text{ Reg } \cap (\mu,\lambda)$
of cardinality  $\le \mu$ we have \newline
${\text{ sup pcf}_{\kappa-\text{complete}}}$
${(\frak a)} \le \lambda$,
\sn
\item  there is no family ${\Cal P}$ of $> \lambda$ subsets of $\lambda$
such that for some regular $\theta \in (\kappa,\mu)$ we have: $A \ne B \in 
{\Cal P} \Rightarrow |A \cap B| < \theta \and |A| \ge \theta$
\sn
\item  cov$(\lambda,\mu^+,\mu^+,\kappa) \le \lambda$
(equivalently cov$(\lambda,\mu,\mu,\kappa) \le \lambda$ as without loss of
generality cf$(\kappa) > \sigma$).
\endroster
Hence
\roster
\item "{(4)}"  there is ${\Cal P} \subseteq [\lambda]^{<\mu}$ such that
$|{\Cal P}| = \lambda$ and every $A \in [\lambda]^{\le \mu}$ is equal to the
union of $< \kappa$ members of ${\Cal P}$
\sn
\item "{(5)}"  there is no tree with $\lambda$ nodes and $> \lambda$ \,
$\theta$-branches when $\theta \in (\kappa,\mu)$ is regular.
\endroster
\enddemo
\bigskip

\demo{Proof}  By \cite{Sh:g}, in detail (we repeat rather than quote
immediate proofs). \nl
1) Without loss of generality cf$(\lambda) \notin [\kappa,\mu)$.  \nl
Note that sup(pcf$_{\kappa\text{-complete}}({\frak a})) \le \sup
\{\text{pp}_{\Gamma(|{\frak a}|^+,\kappa)}(\lambda'):\lambda' = \sup
({\frak a} \cap \lambda')$ and cf$(\lambda') \ge \kappa$ so cf$(\lambda')
\le |{\frak a}| < \mu\}$, and easily the latter is $\le \lambda$ by
\scite{1.1}. \nl
2) By part (4) it is easy (let ${\Cal P}_4 \subseteq [\lambda]^{< \mu}$ be
as in part (4) and $\theta,{\Cal P}_2$ be a counterexample to part (2), so
for every $A \in {\Cal P}_2$ we can find ${\Cal P}'_A \subseteq {\Cal P}_4$
such that $|{\Cal P}'_A| < \kappa$ and $A = \cup\{B:B \in {\Cal P}'_A\}$
hence there is $B_A \in {\Cal P}'_A$ such that $|B_A| = \theta$.  So
$A \mapsto B_A$ is a function from ${\Cal P}_2$ into ${\Cal P}_4$ and
$B_A \in [A]^\theta$ and $A_1 \ne A_2 \in {\Cal P}_2 \Rightarrow |A_1 \cap
A_2| < \theta \and \theta \le |A_1| \and \theta \le |A_2|$ so the function is
one-to-one so $|{\Cal P}_2| \le |{\Cal P}_4| \le \lambda$, contradiction). \nl
3) By \cite[Ch.II,5.4]{Sh:g}. \nl
4) Let ${\Cal P}_0 \subseteq [\lambda]^{< \mu}$ be such that
$|{\Cal P}_0| \le \lambda$ and every $A \subseteq [\lambda]^{\le \mu}$ is
included in the union of $< \kappa$ members of ${\Cal P}_0$ (exists by
part (3)).  Define ${\Cal P} = \{B:\text{for some } A \in {\Cal P}_0,
B \subseteq A\}$ so ${\Cal P} \subseteq [\lambda]^{< \mu}$ and
$|{\Cal P}| \le |{\Cal P}_0| \le |{\Cal P}_0| \cdot \sup\{2^{|A|}:A \in
{\Cal P}_0\} \le \lambda \cdot \mu = \lambda$.
\sn
Now for every $A \in [\lambda]^{\le \mu}$ we can find $\alpha < \kappa$
and $B_i \in {\Cal P}_0$ for $i < \alpha$ such that $A \subseteq
\dbcu_{i < \alpha} B_i$.  Let $B'_i = A \cap B_i$ for $i < \alpha$ so
$B'_i \in {\Cal P}$ and $A = \dbcu_{i < \alpha} B'_i$ as required. \nl
5) Follows by part (2): if the tree is $T$, \wilog \, its set of nodes is
$\subseteq \lambda$ and the set of $\theta$-branches cannot serve as a
counterexample.  \hfill$\square_{\scite{1.2}}$
\enddemo
\bigskip

\remark{\stag{1.2A} Remark}  We can let 
$\mu$ be regular (strong limit $> \aleph_0$)
if we restrict ourselves in \scite{1.2}(1) to 
$|{\frak a}| < \mu$, and in \scite{1.2}(3),(4) to $A \in [\lambda]^{<\mu}$
as if for $\mu' \in \{\mu' < \mu:\mu' \text{ strong limit singular}\},
\kappa(\mu',\lambda)$ is as in \scite{1.2}, then by Fodor's lemma for some
$\kappa = \kappa(\lambda)$ the set 
$S'_\kappa = \{\mu' < \mu:\kappa(\mu',\lambda) =
\kappa\}$ is stationary: this $\kappa$ can serve.
\endremark
\bigskip

\noindent  The stimulation for proving this was in \cite{Sh:454a}
where we actually use:
\demo{\stag{1.3} Conclusion}  Assume $\mu$ is strong limit, $\lambda \ge \mu$.
\ub{Then} for some $\kappa < \mu$ and family ${\Cal P}$,  $|{\Cal P}| \le
\lambda$ we have: for every  $n < \omega$ and $\sigma \in (\kappa,\mu)$
and $f:\left[ \beth_n(\sigma)^+ \right]^{n+1} \rightarrow \lambda$,
for some  $A \subseteq \beth_n(\sigma)^+$ of cardinality $\sigma^+$ we have
$f \restriction A \in {\Cal P}$.
\enddemo
\bigskip

\demo{Proof}  Let $\kappa$ be as in \scite{1.2} (or \scite{1.2A}), and
${\Cal P}$ as in \scite{1.2}(4), and let \newline
${\Cal P}_1 = \{ f:f \text{ a function from some bounded subset } A
\text{ of } \mu \text{ into some } B \in {\Cal P}$ \newline
$\text{(hence } |B| < \mu) \}$. As $\mu$
is strong limit and $|{\Cal P}| \le \lambda$, $\mu \le \lambda$ clearly
$|{\Cal P}_1| \le \lambda$.  Now for any given $f:\left[ \beth_n(\sigma)^+
\right]^{n+1} \rightarrow \lambda$, we can find $\alpha < \kappa$ and
$B_i \in {\Cal P}$ for $i < \alpha$ such that Rang$(f) \subseteq \dsize
\bigcup_{i < \alpha} B_i$.  Define $g:\left[ \beth_n(\sigma)^+ \right]^{n+1}
\rightarrow \alpha$ by: $g(w) = \text{Min} \{ i < \alpha:f(w) \in B_i \}$,
so by the Erd\"os-Rado theorem for some $A \subseteq \beth_n(\sigma)^+$, 
we have: $|A| = \sigma^+$ and $g \restriction A$ is constantly $i(*)$.  Now
$f \restriction A \in {\Cal P}_1$ so we have finished.
\hfill$\square_{\scite{1.3}}$
\enddemo
\bigskip

\demo{\stag{1.4} Conclusion}  If $\lambda = \aleph_0$ or $\lambda$ 
strong limit of cofinality  $\aleph_0$, $(\Omega,{\Cal T})$ is a 
topology (i.e.
$\Omega$ the set of points, ${\Cal T}$ the family of open sets; the topology
is not necessarily Hausdorff or even $T_0)$, ${\Cal B} \subseteq {\Cal T}$
a basis (i.e. every member of ${\Cal T}$ is the union of some subfamily of
${\Cal B}$), and $|{\Cal T}| > |{\Cal B}| + \lambda$ \underbar{then}
$|{\Cal T}| \ge 2^\lambda$.
\enddemo
\bigskip

\demo{Proof}  By \cite{Sh:454a} - the only missing point is that for
$\lambda > \aleph_0$, we need: for arbitrarily large $\mu < \lambda$
there is $\kappa \in (\beth_2(\mu)^+,\lambda)$ such that
cov$(|B|,\kappa^+,\kappa^+,\mu) \le |B|$, which holds by \scite{1.1} (really
in the proof there we use \scite{1.3}). \hfill$\square_{\scite{1.4}}$
\enddemo
\bigskip

\demo{\stag{1.5} Proof of 1.1}  Assume this fails.
By Fodor's Lemma (as in \scite{1.2A}) 
without loss of generality cf$(\mu) = \aleph_0$. 

Without loss of generality for our given $\mu,\lambda$ is the minimal 
counterexample.  Let  $\mu = \dsize \sum_{n<\omega} \mu_n,\mu_n = 
\text{ cf}(\mu_n) < \mu$;
so for each $n$ there is $\lambda_n \in (\mu,\lambda)$ such that
pp$_{\Gamma(\mu^+,\mu_n)}(\lambda_n) \ge \lambda$; hence for some
${\frak a}_n \subseteq \text{ Reg } \cap(\mu,\lambda_n)$ of cardinality
$\le \mu$ and $\mu_n$-complete ideal $J_n \supseteq 
J^{\text{bd}}_{{\frak a}_n}$ we have $\lambda_n = \sup({\frak a}_n)$ and 
$\Pi {\frak a}_n/J_n$ has true cofinality
which is $\ge \lambda$.  Let $\theta_n = \text{ cf}(\lambda_n)$, so $\mu_n \le
\theta_n \le |{\frak a}_n|$.

Without loss of generality $\mu_n > \aleph_0$ hence without loss of 
generality  $|{\frak a}_n| < \mu$ hence without loss of generality
$|{\frak a}_n| < \mu_{n+1}$ (and really even $|\text{pcf}({\frak a}_n)| < 
\mu_{n+1})$, hence the $\theta_n$'s are distinct hence the $\lambda_n$'s 
are distinct, and without loss of generality for $n < \omega$ we have
$\lambda_n < \lambda_{n+1}$ and $\theta_n < \theta_{n+1} < \mu$, hence 
necessarily (by $\lambda$'s minimality)
$\lambda = \dsize \sum_{n<\omega} \lambda_n$, hence without loss of
generality (see \cite[5.2]{Sh:E12}) tcf$\left( \Pi{\frak a}_n,\le_{J_n} 
\right) = \lambda^+$. 

It is clear that forcing by a forcing notion  $Q$  of cardinality  $< \mu$
changes nothing, i.e., we have the same minimal $\lambda$, etc. 
(only omit some $\mu_n$'s).  So without loss of generality $\mu_0 = 
\theta_0 = |{\frak a}_0| = |\text{pcf}({\frak a}_0)| = \aleph_1$,
and for some increasing sequence $\langle \sigma_i:i < \omega_1 \rangle$ of
regular cardinals $< \lambda_0$
\mr
\item "{$(*)$}"  $\lambda_0 = \dsize \sum_{i<\omega_1} \sigma_i$ and
$\dsize \prod_{i<\omega_1} \sigma_i/{\Cal D}_{\omega_1}$
has true cofinality $\lambda^+$ \nl
$({\Cal D}_{\omega_1}$ is the club filter on $\omega_1).$
\ermn
(Of course, we can alternatively use the generalization of 
normal filters as in \cite[\S5]{Sh:410} hence avoid forcing).  
(How do we force?  First by Levy$(\aleph_0,< \mu_0)$ then
Levy$(\mu_0,|\text{pcf}({\frak a}_0)|)$; there is no change in the pcf 
structure for a set of cardinals $> |\text{ pcf}({\frak a}_0)|$, 
so now $|{\frak a}_0| =
\aleph_1$, sup pcf$_{\aleph_1\text{-complete}}({\frak a}_0) > \lambda$ and
pcf$({\frak a}_0)$ has cardinality $\aleph_1$, let ${\frak a}_0 = 
\{\tau_\varepsilon:\varepsilon < \omega_1\}$, pcf$({\frak a}_0) =
\{\theta_\varepsilon:\varepsilon < \omega_1\}$, choose by induction
$\zeta(\varepsilon) < \omega_1$ such that $\tau_{\zeta,\varepsilon} \notin
\cup\{{\frak b}_{\theta_\xi}[{\frak a}_0]:\xi < \varepsilon \text{ and }
\theta_\varepsilon < \lambda\}$, so $\dsize \prod_{\varepsilon < \omega_1}
\theta_{\zeta(\varepsilon)}/J^{\text{bd}}_{\omega_1}$ is $\lambda^+$-directed,
so we get $(*)$ and the statement before it).
Without loss of generality 
\mr
\item "{$(*)_1$}"  $\alpha < \mu_n \Rightarrow |\alpha|^{\aleph_1}
+ \beth_3(\aleph_1) < \mu_n$ for $n \ge 1$.
\ermn
Now by \cite[Ch.VI,\S1]{Sh:g} there is a forcing notion $Q$ of cardinality
$\beth_3(\aleph_1)\,(< \mu !)$ and a name $\underset\tilde {}\to D$ of an
ultrafilter on the Boolean Algebra ${\Cal P}(\omega_1)^V$ (i.e. not
on subsets of $\omega_1$ which forcing by $Q$ adds) which is normal (for
pressing down functions from $V$), extends ${\Cal D}_{\omega_1}$ and, the
main point, the ultrapower $M =: V^{\omega_1}/{\underset\tilde {}\to D}$
(computed in $V^Q$ but the functions are from $V$) satisfies:
\medskip
\roster 
\item "$(*)_2$"  for every  $\kappa > \beth_3(\aleph_1)$ regular or at least
cf$(\kappa) > \beth_3(\aleph_1)$, for some \nl
$g_\kappa \in {}^{\omega_1}$Ord from $V$ (but depending on 
the generic subset of $Q$), the set \newline
$\{ g/ {\approx_{\underset\tilde {}\to D}}:g \in \left(
{}^{\omega_1}\text{Ord}\right)^V, g {<_{\underset\tilde {}\to D}} g_\kappa\}$
is $\kappa$-like (i.e. of cardinality $\kappa$ but every proper
initial segment has cardinality $< \kappa$), the order being
$<_{\underset\tilde {}\to D}$ of course.  We shall say in short
``$g_\kappa / {\underset\tilde {}\to D}$ is $\kappa$-like'', note that for
each  $\kappa$  there is at most one such member in  $M$ (as the ``ordinals''
of $M$ are linearly ordered).
\endroster
\medskip

\noindent
However, we should remember  $V^{\omega_1}/{\underset\tilde {}\to D}$ is,
in general, not well-founded; still there is a canonical elementary 
embedding $\bold j$ of $V$ into $M = V^{\omega_1}/ \underset\tilde {}\to D$
(of course it depends on $G$).
Note that $\bold j$  maps the natural numbers onto  $\{ x \in M:M \models
``x \in \bold j(\omega)"\}$, but this fails for $\omega_1$; \wilog \,
$\bold j \restriction (\omega +1)$ is the identity.  If $M \models
``x$  an ordinal'' let card$_M(x)$  be the cardinality in $V^Q$ of
$\{y:M \models y < x \}$.  Note: also $\bold j(\mu)$ is $\mu$-like and
$\{ \bold j(\mu_n):n < \omega \}$ is unbounded in $\bold j(\mu).$ 

Without loss of generality for every $n \ge 1,\mu_n > |Q|$, and 
$\text{Min}({\frak a}_{n+1}) > \lambda_n$.  For
every regular  $\kappa \in (\mu_1,\lambda^+]$ there is $x_\kappa = g_\kappa /
{\underset\tilde {}\to D}$ which is $\kappa$-like. Note: $g_\kappa \in V$ (not
$\in V^Q \backslash V)$, but we need the generic subset of $Q$ to know which
member of $V$ it is.  Let $\{ g_{\kappa,i}: i < i_\kappa \} \in V$ be a set
such that $\Vdash_Q$ ``for some $i < i_\kappa$ we have $g_{\kappa,i}/
\underset\tilde {}\to D$ is $\kappa$-like" and $i_\kappa \le 
\beth_3(\aleph_1)$.  For regular (in $V$) cardinal 
$\kappa \in (\mu,\lambda^+)$, necessarily  
$M \models ``x_\kappa$ is regular $ > \bold j(\mu)$ and 
$\le g_{\lambda^+}/{\underset\tilde {}\to D}"$ hence without
loss of generality $g_{\lambda^+} = \langle \sigma_\varepsilon:\varepsilon
 < \omega_1 \rangle$ (why? see $(*)$, by \cite[Ch.V]{Sh:g} for some normal
filter ${\Cal D}$ on $\omega_1$ and $\sigma'_\varepsilon \le 
\sigma_\varepsilon$ we have $\dsize \prod_{\varepsilon < \omega_1}
\sigma'_\varepsilon/{\Cal D}$ is $\lambda^+$-like, and force as above; by
renaming we have the above). 
\sn

Now also \wilog \, for regular $\kappa \in (\mu,\lambda^+]$ and 
$i < i_\kappa$ we have Rang$(g_{\kappa,i})$
is a set of regular cardinals $ > \mu$ but $< \lambda_0$
of cardinality $\aleph_1$ (as without loss of generality 
$g_{\kappa,i}(\varepsilon) < \sigma_\varepsilon$ for 
$\varepsilon < \omega_1$ and recall
$\sigma_\varepsilon < \lambda_0$).  For $n \ge 1$ denote 
${\frak c}_n =: \cup \{\text{Rang}(g_{\kappa,i}): \kappa \in {\frak a}_n,
i < i_\kappa\}$ and ${\frak d}_n =: 
\bold j({\frak c}_n) \in M$; note $V \models ``|{\frak c}_n| \le 
|{\frak a}_n| + |Q| = |{\frak a}_n|"$.  So $M \models ``{\frak d}_n$ 
is a set of regular cardinals, each $ > \bold j(\mu)$ but $< \bold j
(\lambda_0)$, of cardinality $\le \bold j(|{\frak a}_n|) < \bold j(\mu_{n+1})
< \bold j(\mu)"$.  Also for every $\kappa \in {\frak a}_n$ we have
$M \models ``x_\kappa \in {\frak d}_n"$ as 
$x_\kappa = g_{\kappa,i} /{\underset\tilde {}\to D}$ for some $i < i_\kappa$
and Rang$(g_{\kappa,i}) \subseteq {\frak c}_n$.
\medskip
\noindent
We can apply the theorem on the structure of pcf (\cite[Ch.VIII,2.6]{Sh:g})
in $M$ (as $M$ is elementarily equivalent to $V$) and get
$\langle {\frak b}_y[{\frak d}_n]:y \in \text{ pcf}({\frak d}_n) \rangle 
\in M$ and $\left< \langle {f^{{\frak d}_n,y}_t}:t < y
\rangle:y \in \text{ pcf}({\frak d}_n) \right> \in M$ 
(this is not a real sequence, only $M$ ``thinks'' so). 
\medskip
For  $y \in M$  such that $M \models ``y$  a limit ordinal (e.g. a
cardinal)'' let  $\lambda_y$ be the cofinality (in $V^Q)$ of $( \{ x:M
\models ``x \text{ an ordinal } < y" \},<^M)$. So
\mr
\item "{$(*)_3$}"  $\kappa = \lambda_{(x_\kappa)}$ for $\kappa \in
\text{ Reg}, \kappa > |Q|$ 
\sn
\item "{$(*)_4$}"  assume $|\{a:a \in^M \bold j(\mu_m)\}| < \mu_n$, then \nl
$M \models ``\sup \text{ pcf}_{\bold j(\mu_m)\text{-complete}}
({\frak d}_n \cap g_{\lambda^+}/ \underset\tilde {}\to D) \ge g_{\lambda^+}
/ \underset\tilde {}\to D$". 
\ermn
[Why?   Assume not, so $M \models ``\text{sup pcf}_{\bold j(\mu_m)
\text{-complete}}({\frak d}_n \cap g_{\lambda^+}/\underset\tilde {}\to D) <
g_{\lambda^+}/\underset\tilde {}\to D$" hence $M \models$ ``for every
$g \in \Pi({\frak d}_n \cap g_{\lambda^+}/\underset\tilde {}\to D)$ for some 
$\langle (y_\ell,a_\ell):\ell < \bold j(\mu_m) \rangle,y_\ell \in 
\text{ pcf}({\frak d}_n \cap g_{\lambda^+}/\underset\tilde {}\to D),a_\ell$ 
an ordinal $< y_\ell$ we have 
$g < \underset {\ell < \bold j(\mu_m)}\to {\text{sup}} 
f^{y_\ell}_{a_\ell}$".  In $V^Q$ we have $\Pi {\frak a}_n/J_n$ is 
$\lambda^+$-directed hence $\dsize \prod_{\kappa \in {\frak a}_n} (\{t:t <^M 
x_\kappa\},<^M)/J_n$ is $\lambda^+$-directed (by $(*)_3$) hence there is a 
function $g^*$ such that
\mr
\item "{$(a)$}"  Dom$(g^*) = {\frak a}_n$
\sn
\item "{$(b)$}"  $g^*(\kappa) <^M x_\kappa = 
g_\kappa/\underset\tilde {}\to D$
\sn
\item "{$(c)$}"  if $M \models ``y \in \text{ pcf}({\frak d}_n \cap 
g_{\lambda^+}/\underset\tilde {}\to D)$ and $a < y"$ then \nl
$\{\kappa \in {\frak a}_n:M \models ``f^{{\frak d}_n,y}_a(x_\kappa) <^M 
g^*(\kappa)"\} = {\frak a}_n \text{ mod } J_n$.
\ermn
By \scite{1.6}(2) below we can find $Y \in V$ such that $|Y| < |Q|^+ + \mu 
= \mu$ and $\kappa \in {\frak a}_n \Rightarrow M \models ``g^*(\kappa) \in 
\bold j(Y)"$.
There is $g^\otimes \in M$ such that $M \models ``g^\otimes \in \Pi
{\frak d}_n$ and $g^\otimes(\theta) = (\sup(\bold j(Y)) \cap \theta) + 1
< \theta"$ (as $M \models ``\text{Min}({\frak d}_n) > \bold j(\mu)"$). \nl
By the choice of $Y$ clearly $\kappa \in {\frak a}_n \Rightarrow g^*(\kappa)
<^M g^\otimes(\kappa)$. \nl
By the choice of $\left< \langle f^{{\frak d}_n,y}_t:t < y \rangle:y \in
\text{ pcf}({\frak d}_n) \right>$ (in $M$'s sense) and the assumption toward
contradiction we have: \nl

$M \models$ ``there is a subset $\Theta$ of pcf$({\frak d}_n \cap
g_{\lambda^+}/\underset\tilde {}\to D)$ of cardinality 
$< \bold j(\mu_m)$ and $\langle a_\theta:
\theta \in \Theta \rangle \in \Pi \Theta$ such that $(\forall \sigma \in
{\frak d}_n)(\dsize \bigvee_{\theta \in \Theta} g^\otimes(\sigma) <
f^{{\frak d}_n,\theta}_{a_\theta}(\sigma))$".  Choose such a sequence $\langle
a_\theta:\theta \in \Theta \rangle$ in $M$ and let $\langle \theta_i:i < i
(*) \rangle$ list the $\theta \in^M \Theta$, so $i(*) < \mu_n$ by the
hypothesis of $(*)_4$.  Let
${\frak a}_{n,i} = \{\kappa \in {\frak a}_n$:letting $\sigma = x_\kappa \in
\mu$ we have $g^*(\sigma) < f^{{\frak d}_n,\theta_i}_{a_{\theta_i}}(\sigma)\}
\in V^Q$.  Now as $g^*(\kappa) < g^\otimes(x_\kappa)$, clearly ${\frak a}_n
= \dbcu_{i < i(*)} {\frak a}_{n,i}$.  So for some $i < i(*)$ we have 
${\frak a}_{n,i} \in J^+_n$, and we get a contradiction to
the choice of $g^*$ hence at last we have proved $(*)_4$.]
\sn
Clearly $\bold j(\langle {\frak c}_n:n < \omega \rangle )$ is a sequence
of length $\bold j(\omega) = \omega$ hence $\bold j(\langle {\frak c}_n:n <
\omega \rangle) = \langle {\frak d}_n:n < \omega \rangle$, i.e. with $n$-th
element ${\frak d}_n$.
Let ${\bar z} \in M$ be such that $M \models ``{\frak p} = \langle
(k_n,t_n,s_n):n < \omega \rangle$ defined by: $k_n < \omega$ is maximal such
that $g_{\lambda^+}/\underset\tilde {}\to D \le \text{ sup pcf}
_{\bold j(\mu_n)\text{-complete}}({\frak d}_n \cap g_{\lambda^+}/
\underset\tilde {}\to D)$, and $t_n$ is the minimal cardinal such
that sup pcf$_{\bold j(\mu_n)\text{-complete}}$ (${\frak d}_n \cap 
(g_{\lambda^+}/\underset\tilde {}\to D))$ is 
$\ge g_{\lambda^+}/\underset\tilde {}\to D$ and cf$(t_n) = s_n$ so
$s_n \ge \bold j(\mu_n)"$.  As $\bold j(\mu)$ is $\mu$-like
clearly $(\forall m < \omega)(\exists n < \omega)(m < n \and |\{x \in M:
x \in^M (\bold j(\mu_m))\}| < \mu_n)$ hence by $(*)_4$ above
necessarily $(\forall m < \omega)(\exists n < \omega)
\left[ |[s_n]| \ge \mu_m \right]$, but $\bold j(\mu)$ is
the limit of $\langle \bold j(\mu_n):n < \omega \rangle \in M$,
hence $M \models ``\bold j(\mu) = \lim z_n"$.  Now 
\mr
\item "{$(*)_5$}"  $M \models ``\bold j(\mu),g_{\lambda^+}/
\underset\tilde {}\to D$ form a counterexample to the Theorem \scite{1.1}".
\ermn
But as $\bold j$ is an elementary embedding of $V$ to $M$, the choice 
of $\lambda$ (minimal) implies
$$
\align
M \models &``\text{ there is no } \lambda' < \bold j(\lambda)
 \text{ such that }
\bold j(\mu),\lambda' \\
          &\text{ form a counterexample to the theorem}".
\endalign
$$
\noindent
But as $\text{ Rang }[g_{\lambda^+}/{\underset\tilde {}\to D}] < \bold j
(\mu_0)< \bold j(\lambda)$, clearly we have $M \models 
``g_\lambda / {\underset\tilde {}\to D} < \bold j(\lambda)"$.
\newline
\noindent
By the last two sentences we get a contradiction to $(*)_5$. 
\hfill$\square_{\scite{1.1}}$
\enddemo
\bigskip

\demo{\stag{1.6} Observation}  Let 
$Q,{\underset\tilde {}\to D},G \subseteq Q,V^Q,M,\bold j$ 
be as in the proof \scite{1.5}.  
Let for $z \in M$,  $[z] = \{ t:M \models t \in y \}$.  So
\roster
\item  \underbar{If}  $Y \in V^Q$,  $Y \subseteq M$, $\chi = \text{Max}
\left\{ |Y|^{V^Q},|Q|^V \right\}$ \underbar{then} for some $y \in V$,
$|y|^V \le \chi$ and $\forall x[x \in Y
\Rightarrow  M \models ``x \in \bold j(y)"]$.
\sn
\item  Assume  $M \models ``{\frak d}$ is a set of regular cardinals
$ > |{\frak d}|, > \bold j \left( |Q|^V \right)"$ and $\lambda_y$ (when
$M \models$ ``$y$ limit ordinal'') is as in \scite{1.5} 
(its cofinality in  $V^Q$).
{\roster
\itemitem{ (a) } If  $M \models ``y \in \text{ pcf}({\frak d})",J$ is 
(in $V^Q$) the ideal on $[{\frak d}]$ generated by \nl
$\{ [{\frak b}_\theta [{\frak d}]]:M \models
``\theta \in \text{ pcf}({\frak d}) \and \theta < y"\} \cup
\{[{\frak d} \setminus {\frak b}_y[{\frak d}]] \}$ \nl
\underbar{then} (in $V^Q$) 
$\dsize \prod_{x \in [{\frak d}]} \lambda_x/J$  has true cofinality
$\lambda_y$
\sn
\itemitem{ (b) }  cf$\left( \Pi \{ \lambda_y:y \in [{\frak d}] \} \right) =
\text{max}\{ \lambda_y:y \in [\text{pcf }{\frak d}] \}$.
\endroster}
\endroster
\enddemo

\demo{Proof}  Straightforward (and we use only part (1)).  
For (2)(b) remember
$$
M \models ``y \text{ is finite }" \Rightarrow [y] \text{ finite}.
$$
\enddemo
\medskip

\remark{\stag{1.7} Remark}  Of course, the proof of \scite{1.1} gives
somewhat more than stated (say after fixing $\mu_0 = \aleph_1$).  E.g.,
\mr
\item "{$\oplus$}"  the cardinal $\mu$ satisfies the conclusion of \scite{1.1}
for $\lambda \ge \lambda^*$ if
\sn
\item "{$\boxtimes_\mu$}"  $\mu > \text{ cf}(\mu) = \aleph_0$ (as before this
suffices) and $\mu = \sup\{\kappa < \mu:\kappa \text{ is regular}$ \nl
uncountable and there is a forcing notion $Q \text{ satisfying the } 
\kappa\text{-c.c. of cardinality}$ \nl
$\le \lambda_0 < \mu\}$ such that $\Vdash_Q$ ``for every
$\aleph_1$-complete filter $D$ on $\kappa$ from $V$ containing the 
co-countable sets there is an ultrafilter $\underset\tilde {}\to D$ on 
${\Cal P}(\kappa)^V$ extending $D$ as in \cite[Ch.VI,\S1]{Sh:g} for regular
cardinal $> \lambda^+$ which is complete for partitions of $\kappa$ from $V$
to countably many parts.
\ermn
\endremark
Alternatively, we can phrase the theorem after fixing $D$.
\newpage

\head{\S2 The Main Theorem Revisited} \endhead \resetall
\bigskip

We give another proof and get more refined information.  Note that in 
\scite{2.1}
if $\mu$ is strong limit, we can choose $R^\ast$ such that: if $\theta <
\kappa$ are in $R^\ast$ then $2^\theta < \kappa$ and then $\otimes^0_{R^\ast,
\theta,\theta_1}$ is immediate.
\proclaim{\stag{2.1} Theorem}  Suppose  $\mu$  is a limit singular cardinal
satisfying:
\roster
\item "{$\otimes^0_\mu$}"  for any  $R \subseteq \mu \cap \text{\rm Reg}$
unbounded, for some $\theta \in R,\theta > \text{\rm cf}(\mu)$ and $\theta_1$,
{\rm cf}$(\mu) < \theta_1,\aleph_1 \le \theta_1 < \theta$ and $R^* \subseteq
R \backslash \theta^+$ unbounded in $\mu$ we have: 
\endroster
\medskip

\roster
\item "{$\otimes^0_{R^\ast,\theta,\theta_1}$}"
\underbar{if} $\sigma < \kappa$
are in $R^\ast$, $f_\alpha:\theta \rightarrow \sigma$ for $\alpha < \kappa$,
$I_\kappa$ a $\kappa$-complete ideal on $\kappa$ extending 
$J^{\text{\rm bd}}_\kappa$
and $J$ is a $\theta$-complete ideal on $\theta$, \underbar{then} for some 
$A \in I^+_\kappa$ and \nl
$B_\alpha \subseteq \theta$ for $\alpha \in A$ satisfying $B_\alpha = 
\theta \text{ mod }J$ we have \nl
$\xi < \theta \Rightarrow | \{ f_\alpha(\xi):\alpha \in A \text{ and } 
\xi \in B_\alpha\}| < \theta_1$.
\endroster
\medskip

\noindent
\underbar{Then}
$$
\text{ for every } \lambda > \mu \text{ we have: } \tag"{$\otimes^1_\mu$}"
$$
\medskip

$$
\text{ for some } \kappa < \mu \text{ we have: }
\tag"{$\otimes^1_{\lambda,\mu}$}"
$$
\medskip

$$
\text{ for every } {\frak a} \subseteq
(\mu,\lambda) \cap \text{ {\rm Reg\/} of cardinality } < \mu,
\text{{\rm pcf\/}}_{\kappa-{\text{\rm complete\/}}}({\frak a}) 
\subseteq \lambda. \tag"${\otimes^1_{\lambda,\mu,\kappa}}$"
$$
\endproclaim
\bigskip
\noindent
Before we prove it, note:

\demo{\stag{2.2} Observation}  Assume: 
\roster
\item "{(a)}"  $\langle w^n_i:i < \alpha^\ast \rangle$ is a sequence
of pairwise disjoint sets,  $w^n = \dsize \bigcup_{i<\alpha^*}w^n_i$ \nl
(possibly  $w^n_i = \emptyset$  for some  $n$  and $i$)
\sn
\item "{(b)}"  $\left( \underset {n,i}\to {\text{sup}}
|w^n_i|^+ \right) < \theta$ and $\theta$ is uncountable
\sn
\item "{(c)}"  $J_n$ is a $\theta$-complete ideal on $w^n$ such that
$w^n \notin J_n$
\sn
\item "{(d)}"  $h^n_i:w^{n+1}_i \rightarrow w^n_i$ and $h^n = \dsize \bigcup
_{i<\alpha^\ast}h^n_i$
\item "{(e)}"  for every $A \in  J_{n+1}$ the set
$ \{ x \in w^n:(\forall y \in  
w^{n+1})[h^n(y) = x \Rightarrow  y \in  A] \}$  belongs to  $J_n$.
\endroster
\medskip 

\noindent
\underbar{Then} for some  $i$  there are  $x_n \in w^n_i$ such that
$\dsize \bigwedge_n h^n(x_{n+1}) = x_n$.
\enddemo
\bigskip

\remark{\stag{2.2A} Remark}  Hence for the $J_m$-majority of
$y \in w^m$ there is  
$\langle x_n:n < \omega \rangle$  as above such that $y = x_m$.
\endremark
\bigskip

\demo{Proof}  Without loss of generality  $\langle w^n_i:n < \omega ,
i < \alpha^\ast \rangle$  are pairwise disjoint.  Now we define by
induction on the ordinal $\zeta \le \theta$  for each  $i < \alpha^*$ 
a set $u^\zeta_i \subseteq w_i =: \dsize \bigcup_{n < \omega} w^n_i$ by: 
$$
u^\zeta_i = \biggl\{ x \in w_i: x \in  
\dsize \bigcup_{\xi < \zeta} u^\xi_i \text{ or }
(\forall y \in w^{n+1}_i)[h^n_i(y) = x \Rightarrow  y \in \dsize
\bigcup_{\xi < \zeta} u^\xi_i \biggr\}.
$$
\medskip
\noindent
So $\langle u^\zeta_i:\zeta < \theta \rangle$ is an increasing 
sequence of subsets of $w_i$.  Also $u^{\zeta+1}_i = u^\zeta_i \Rightarrow
(\forall \xi > \zeta)[u^\xi_i = u^\zeta_i]$, hence there is for each
$i < \alpha^\ast$ a unique $\zeta[i] < \aleph_1 + |w_i|^+$ such that
$u^\zeta_i =
u^{\zeta[i]}_i \Leftrightarrow \zeta \ge \zeta[i]$. \newline
\noindent
If for some $i$ we have $u^{\zeta[i]}_i \ne w_i$, we can easily prove the 
conclusion so assume $u^{\zeta[i]}_i = w_i$ for every $i$.
Let $\mu = \underset {i}\to {\text{sup}}(|w_i|^+ + \aleph_1)$,
so except when $\theta \le \aleph_1$ (hence $\theta = \aleph_1$) 
we know $\mu < \theta$.  Now we can use clause (e) to prove by induction on  
$\zeta \le \mu$  for all  $n$  that
$$
\cup \{ u^\zeta_i \cap w^n_i:i < \alpha^\ast \} \in J_n
$$
\medskip
\noindent
(we use  $J_n$ is $\theta$-complete, $\theta > \mu)$.  But as 
$i = \mu \Rightarrow u^\mu_i \cap w^n_i = w^n_i$ we get $w^n \in J_n$, 
a contradiction.  
We are left with the case $\theta = \aleph_1$
so each $w^n_i$ is finite and $i < \alpha^* \Rightarrow \zeta[i] < \theta$ 
but then for each $m$ we have $\cup \{ u^m_i \cap w^0_i:i < \alpha^*\} \in  
J_0$, so as $J_0$ is $\theta$-complete
there is  $x \in w^0$ such that for each $m < \omega$ and $i < \alpha^*$ 
we have $x \notin u^m_i \cap w^0_i$.
For some  $i(\ast)$, $x \in  w^0_{i(\ast)}$, so as $x \notin u^n_{i(\ast )}$ 
for some  $x_n \in w^n_{i(\ast)}$ we have $h^{n-1} \circ h^{n-2} \circ  
\cdot \cdot \cdot \circ  h_0(x_n) = x$.  By K\"onig's Lemma (as all  
$w^n_{i(\ast )}$ are finite) we finish. \nl
${{}}$  \hfill$\square_{\scite{2.2}}$
\enddemo
\bn
Before we continue we mention some things which are essentially from
\cite{Sh:g}, and more explicitly, \cite[6.7A]{Sh:430}. \nl
We forgot there to mention the most obvious demand
\proclaim{\stag{2.3A} Subclaim}  In {\rm \cite[6.7A\/]{Sh:430}} we can add:
\mr
\item "{$(j)$}"  ${\text{\rm max pcf\/}}({\frak b}^\beta_\lambda
[\bar{\frak a}]) = \lambda$ (when defined).
\ermn
Also in {\rm\cite[6.7\/]{Sh:430}} we can add
\mr
\item "{$(\delta)$}"  ${\text{\rm max pcf\/}}({\frak b}_\lambda) = \lambda$.
\endroster
\endproclaim
\bigskip

\demo{Proof}  This is proved during the proof of \cite[6.7]{Sh:430} (see
$(*)_4$ in that proof, p.103).  Actually we have to state it earlier in
$(*)_2$ there, i.e. add
\mr
\item "{$(\zeta)$}"  max pcf$({\frak b}^{i,j}_\lambda) \le \lambda$.
\ermn
We then quote \cite[Ch.VIII,1.3,p.316]{Sh:g}, but there this is stated.

Lastly, concerning \cite[6.7A]{Sh:430} the addition is inherited from
\cite[6.7]{Sh:430}.  \nl
${{}}$  \hfill$\square_{\scite{2.3A}}$
\enddemo
\bigskip

\proclaim{\stag{2.3B} Subclaim}  In {\rm\cite[6.7A\/]{Sh:430}} we can deduce:
\mr
\item "{$(\alpha)$}"  if ${\frak a}' \subseteq \dbcu_{\beta < \sigma}
{\frak a}_\beta,|{\frak a}'| < \sigma,{\frak a}' \in N_\sigma$, \ub{then}
for some $\beta(*) < \sigma$ and finite \nl
${\frak c} = \{\theta_0,\theta_1,
\dotsc,\theta_n\} \subseteq {\frak a}_{\beta(*)}$ we have
{\roster
\itemitem{ (i) }  $\theta_\ell > \theta_{\ell +1}$
\sn
\itemitem{ (ii) }  ${\frak a}' \subseteq \dbcu_{\ell \le n}
{\frak b}^{\beta(*)}_{\theta_\ell}[\bar{\frak a}]$
\sn
\itemitem{ (iii) } $\beta \in (\beta(*),\sigma) \Rightarrow {\frak b}^\beta
_{\theta_\ell}[\bar{\frak a}] \cap {\frak a}' = {\frak b}^{\beta(*)}
_{\theta_\ell}[\bar{\frak a}] \cap {\frak a}'$
\sn
\itemitem{ (iv) }  $\theta_\ell = \text{{\rm max pcf\/}}({\frak a}' \backslash
\dbcu_{k < \ell} {\frak b}^{\beta(*)}_{\theta_k}[{\frak a}])$
\endroster}
\item "{$(\beta)$}"  moreover, $\langle \theta_\ell:\ell \le n \rangle$ is
definable from ${\frak a}',\beta(*)$ and $\langle {\frak b}^{\beta(*)}
_\theta[\bar{\frak a}]:\theta \in {\frak a}_\beta \rangle$ uniformly
\sn
\item "{$(\gamma)$}"  if $\langle {\frak a}'_\varepsilon:\varepsilon <
\zeta \rangle \in N_\sigma,\zeta < \sigma,|{\frak a}'_\varepsilon| < \sigma$
\ub{then} we can have one $\beta(*)$ for all ${\frak a}'_\varepsilon$ and so
$\left< \langle \theta_{\varepsilon,\ell}:\ell \le n(\varepsilon) \rangle:
\varepsilon < \zeta \right> \in N_{\beta(*)}$.
\endroster
\endproclaim
\bigskip

\demo{Proof}  \ub{Clause $(\alpha)$}.  We choose $\beta_\ell,\theta_\ell$
by induction on $\ell$.  For $\ell = 0$ clearly for some $\theta_0,
{\frak a}' \in N_{\gamma_0}$ so ${\frak a}' \subseteq {\frak a}_{\gamma_0}$,
hence $\theta_0 = \text{ max pcf}({\frak a}')$ belongs to $N_{\gamma_0}$
hence to ${\frak a}_\beta$ for $\beta \in [\gamma_0,\sigma)$ so by clause
$(i)$ of \cite[6.7A]{Sh:430}, $\langle {\frak b}^\beta_{\theta_0}
[\bar{\frak a}]:\beta \in [\gamma_0,\sigma) \rangle$ is increasing hence
$\langle {\frak b}^\beta_{\theta_0}[\bar{\frak a}] \cap {\frak a}':\beta \in
[\gamma_0,\sigma) \rangle$ is eventually constant, say for $\beta \in
[\beta_0,\sigma),\beta_0 \in (\gamma_0,\sigma)$.  For $\ell +1$ apply the
case $\ell = 0$ to ${\frak a}' \backslash \dbcu_{k \le \ell} 
{\frak b}^{\beta_k}_{\theta_k}[\bar{\frak a}]$ and get $\theta_{\ell +1},
\beta_{\ell +1}$. \hfill$\square_{\scite{2.3B}}$
\sn
\ub{Clauses $(\beta), (\gamma)$}.  Easier.
\enddemo
\bigskip

\proclaim{\stag{2.3} Claim}  1) Assume $\sigma \ge \aleph_0$ is regular, 
$\lambda$ a cardinal, $J$ the $\sigma$-complete ideal generated by 
$J_{<\lambda}[{\frak a}]$ for a set ${\frak a}$ of regular cardinals
$>|{\frak a}|$, ${\frak a} \notin J$, ${\frak a}_i \in J$ for $i < \alpha$,
${\frak a} = \dsize \bigcup_{i<\alpha}{\frak a}_i$ and $\text{{\rm max pcf\/}}
({\frak a}_i) < \lambda$. \nl
\underbar{Then} \footnote{note that \wilog \, $i < \alpha \Rightarrow
{\frak a}_i \ne \emptyset$ so necessarily $|\alpha| \le |{\frak a}|$.} 
we can find ${\frak b}$, ${\frak b}_i(i<\alpha)$ and
$I$ such that:
\roster
\item "{(a)}"  ${\frak b}_i \subseteq \text{{\rm pcf\/}}
({\frak a}_i)$ is finite
\sn
\item "{(b)}"  ${\frak b} = \dsize \bigcup_{i<\alpha}{\frak b}_i$
\sn
\item "{(c)}"  $I$ is an ideal on  ${\frak b}$ 
\sn
\item "{(d)}"  for $w \subseteq \alpha$ we have $\dsize \bigcup_{i \in w}
{\frak a}_i \in J \Leftrightarrow \dsize \bigcup_{i \in w}{\frak b}_i \in I$
\sn
\item "{(e)}"  $I$ is the $\sigma$-complete ideal generated by
$J_{<\lambda}[{\frak b}]$
\sn
\item "{(f)}"  we have ${\frak b}_i = \{ \lambda_{i,\ell}:\ell < n_i \}$ and
if $I_1$ is an $\aleph_1$-complete ideal on ${\frak b}$ extending $I$ (so
$I_1 = I$ is O.K. if $\sigma > \aleph_0)$, \ub{then} for any 
${\frak d} \in I^+_1$ 
there are $B \subseteq \alpha$ and $\ell^\ast < \omega$ such that:
{\roster
\itemitem{ $(\alpha)$ }  $\{ \lambda_{i,\ell^\ast}:i \in B \}
\subseteq {\frak d}$,
\sn
\itemitem{ $(\beta)$ }  $\{ \lambda_{i,\ell^\ast}:i \in B \} \in I^+_1$
\sn
\itemitem{ $(\gamma)$ } $\text{ for every } B^\prime \subseteq B$ we have
$\dsize \bigcup_{i \in B^\prime} {\frak b}_i \in I_1 \Leftrightarrow
\left\{ \lambda_{i,\ell^\ast}:i \in  B^\prime \right\} \in I_1$.
\endroster}
\endroster

\noindent
2)  Assume in addition $\text{{\rm pcf\/}}_{\kappa_i-\text{{\rm complete\/}}}
({\frak a}_i) \subseteq \lambda_i$ and $\kappa_i \le \sigma$ \underbar{then} 
we can find ${\frak b}$, \nl
${\frak b}_i(i < \alpha)$ and $I$ such that:
\roster
\item "{$(a)'$}"  ${\frak b}_i \subseteq 
\text{{\rm pcf\/}}_{\kappa_i-\text{complete}}
({\frak a}_i) \subseteq \lambda_i$ has cardinality $< \kappa_i$
\endroster
and (b) -- (e) hold.
\mn
3)  Assume
\mr
\widestnumber\item{$(iii)$}
\item "{$(i)$}"  $I$ an ideal on $\alpha$
\sn
\item "{$(ii)$}"  $J$ an ideal on $\beta$
\sn
\item "{$(iii)$}"  $\langle \chi_i:i < \alpha \rangle$ a sequence of regular
cardinals with tcf$(\dsize \prod_{i < \alpha} \chi_i/I) = \chi$
\sn
\item "{$(iv)$}"  for $i < \alpha,\langle \tau^i_j:j < \beta \rangle$ is
a sequence of regular cardinals with \nl
tcf$(\dsize \prod_{j < \beta}\tau^i_j/J) = \chi_i$
\sn
\item "{$(v)$}"  $\langle \sigma_j:j < \beta \rangle$ is a sequence of
regular cardinals
\sn
\item "{$(vi)$}"  $|\alpha| + |\beta| + \dsize \sum_{j < \beta} \sigma_j
< \text{ min}\{\tau^i_j:i < \alpha,j < \beta\}$.
\ermn
\ub{Then} there are for each $j < \beta$ an ordinal 
$\varepsilon_j < \sigma_j$ and
sets $\langle {\frak b}^j_\varepsilon:\varepsilon < \varepsilon_j \rangle$
such that
\mr
\item "{$(a)$}"  $\dbcu_{\varepsilon < \varepsilon_j} {\frak b}^j
_\varepsilon \subseteq \{\tau^i_j:i < \alpha \}$ and if max pcf$\{\tau^i_j:
i < \alpha,j < \beta\} = \chi$ then equality holds
\sn
\item "{$(b)$}"  $\lambda^j_\varepsilon =: \text{{\rm max pcf\/}}
({\frak b}^j_\varepsilon)$ is in 
$\text{{\rm pcf\/}}_{\sigma_j\text{-complete}}({\frak b}^j_\varepsilon)$, 
and letting $J^*$ be the ideal with domain 
$\dbcu_{j < \beta} \{j\} \times \varepsilon_j$
defined by $A \in J$ iff $\text{{\rm max pcf\/}}\{\lambda^j_\varepsilon:
(j,\varepsilon) \in A\} < \chi$, we have
\sn
\item "{$(c)$}"  if $w \in J^*$, then $\{i < \alpha:\{j < \beta:\exists
\varepsilon < \varepsilon_x[\tau^i_j \in {\frak b}^j_\varepsilon \wedge
(j,\varepsilon) \in w\} \notin J\} \in I$.
\ermn
(Note that $J^*$ is a proper ideal and $\dsize \prod_{(j,\varepsilon) \in
\text{ Dom}(J^*)} \lambda^j_\varepsilon/J^*$ is $\chi$-directed by basic
$\text{{\rm pcf\/}}$ theory.)
\endproclaim
\bigskip 

\demo{Proof}  By the proof of \cite[Ch.VIII,1.5]{Sh:g} or by 
\cite[6.7,6.7A,6.7B]{Sh:430} (for (1)(f), shrink $A$ to make 
$n_i$ constantly $n^*$, then prove by induction on $n^*$).  In more detail:
\nl
1) Without loss of generality Min$({\frak a}) > |{\frak a}|^{+3}$.
To be able to use \cite{Sh:430} freely in its notation rename ${\frak a}_i$
as ${\frak e}_i$.   We apply \cite[6.7A,p.104]{Sh:430} with 
${\frak a},\kappa,\sigma$ there standing
for ${\frak a},|{\frak a}|^{++},|{\frak a}|^+$ here and \wilog \,
$\langle {\frak e}_i:i < \alpha \rangle \in N_0,\lambda \in N_0$.  By the
subclaim \scite{2.3B} above for each $i < \alpha$ there are $\beta(i) <
|{\frak a}|^+$ and finite ${\frak b}_i \subseteq \text{ pcf}({\frak e}_i)
\cap {\frak a}_{\beta(i)}$ such that $\beta \in [\beta(i),|{\frak a}|^+)
\Rightarrow {\frak e}_i \subseteq \dbcu_{\mu \in {\frak b}_i}
{\frak b}^{\beta +1}_\mu[\bar{\frak a}]$.  Moreover $\langle ({\frak b}_i,
\beta(i)):i < \alpha \rangle \in N_{\beta(*)}$ for $\beta(*) =
(\underset {i < \alpha}\to \sup \beta(i)) +1$ and let ${\frak b} = 
\dbcu_{i < \alpha} {\frak b}_i$ and $J = \{{\frak c} \subseteq {\frak b}:
\text{we can find } \zeta < \sigma \text{ and } \langle {\frak c}_\varepsilon:
\varepsilon < \zeta \rangle$ such that ${\frak c} = \dbcu_{\varepsilon <
\zeta} {\frak c}_\varepsilon$ and max pcf$({\frak c}_\varepsilon) 
< \lambda\}$.  Let us check all the clauses of the desired conclusion.
\mn
\ub{Clause (a)}:  ${\frak b}_i \subseteq \text{ pcf}({\frak e}_i)$ is finite.
\nl
Holds by the choice of ${\frak b}_i$.
\mn
\ub{Clause (b)}:  ${\frak b} = \dbcu_{i < \alpha} {\frak b}_i$. \nl
Holds by the choice of ${\frak b}$.
\mn
\ub{Clause (c)}:  $J$ an ideal on ${\frak b}$. \nl
By \cite[Ch.I]{Sh:g} and the definition of $J$.
\mn
\ub{Clause (d)}:  For $w \subseteq \alpha$ we have $\dbcu_{i \in w}
{\frak e}_i \in J \Leftrightarrow \dbcu_{i \in w} {\frak b}_i \in I$. \nl
Why?  By the definition of $J$ and it suffices to prove for each subset
$w$ of $\alpha$ that

$$
\text{max pcf}(\dbcu_{i \in w} {\frak e}_i) < \lambda \Leftrightarrow
\text{ max pcf}(\dbcu_{i \in w} {\frak b}_i) < \lambda.
$$
\mn
First assume max pcf$(\dbcu_{i \in w} {\frak e}_i) < \lambda$.  Now
$j \in w \Rightarrow {\frak b}_j \subseteq \dbcu_{i \in w} {\frak e}_i$
hence (by \cite[Ch.I,1.11]{Sh:g}) pcf$(\dbcu_{j \in w} {\frak b}_j)
\subseteq \text{ pcf}(\dbcu_{i \in w} {\frak e}_i)$ so
max pcf$(\dbcu_{i \in w} {\frak b}_i) \le \text{ max pcf}(\dbcu_{i \in w}
{\frak e}_i) < \lambda$, as required.
\sn

If the other implication fails, then there is $w \subseteq \alpha$ which
exemplifies it in $N_{\beta(*)}$ (as all the relevant parameters are in it),
so we need only consider $w \in N_{\beta(*)}$.  Assuming $w \in
N_{\beta(*)}$ and max pcf$(\dbcu_{i \in w} {\frak b}_i) < \lambda$ 
let ${\frak b}' =:
\dbcu_{i \in w} {\frak b}_i$, so ${\frak b}' \in N_{\beta(*)} \cap
{\frak a}_{\beta(*)}$ and max pcf$({\frak b}') < \lambda$, and by
\cite[6.7A,(h)]{Sh:430} for some finite ${\frak c} \subseteq \text{ pcf}
({\frak b}') \cap N_{\beta(*)}$ we have $\dbcu_{\theta \in {\frak c}}
{\frak b}^{\beta(*)}_\theta[\bar{\frak a}]$ includes ${\frak b}'$. \nl
By \cite[6.7A(f)]{Sh:430}, i.e., smoothness

$$
\tau \in {\frak b}' \Rightarrow {\frak b}^{\beta(*)}_\tau[\bar{\frak a}]
\subseteq \dbcu_{\theta \in {\frak c}} {\frak b}^{\beta(*)}_\theta
[\bar{\frak a}]
$$
\mn
hence

$$
\align
\tau^* \in \dbcu_{i \in w} {\frak e}_i &\Rightarrow \dsize \bigvee_{i \in w}
\tau^* \in {\frak e}_i \\
  &\Rightarrow \dsize \bigvee_{i \in w} \tau^* \in \cup\{{\frak b}^{\beta(*)}
_\tau[\bar{\frak a}]:\tau \in {\frak b}_i\} \Rightarrow \dsize \bigvee
_{i \in w} \tau^* \in \dbcu_{\theta \in {\frak c}} {\frak b}^{\beta(*)}
_\theta[\bar{\frak a}] \\
  &\Rightarrow \tau^* \in \dbcu_{\theta \in {\frak c}}
{\frak b}^{\beta(*)}_\theta[\bar{\frak a}].
\endalign
$$
\mn
So $\dbcu_{i \in w} {\frak a}_i \subseteq \dbcu_{\theta \in {\frak c}}
{\frak b}^{\beta(*)}_\theta[\bar{\frak a}]$ hence

$$
\align
\text{max pcf}(\dbcu_{i \in w} {\frak e}_i) &\le \text{ max pcf}
(\dbcu_{\theta \in {\frak c}} {\frak b}^{\beta(*)}_\theta[\bar{\frak a}]) \\
  &\le \underset {\theta \in {\frak c}}\to {\text{max}}(\text{max pcf }
{\frak b}^{\beta(*)}_\theta[{\frak a}]) = \text{ max}({\frak c}) < \lambda
\endalign
$$
\mn
(we use subclaim \scite{2.3B} above).
\enddemo
\bn
\ub{Clause (e)}:  $I$ is the $\sigma$-complete ideal generated by
$J_{< \lambda}[{\frak b}]$. \nl
By the choice of $I$.
\bn
\ub{Clause (f)}:  As $I_1$ is $\aleph_1$-complete for some $n^*$ the set
${\frak d} \cap \cup\{{\frak b}_i:|{\frak b}_i| = n^*\}$ belongs to
$I^+_1$.  Now we try to choose by induction on $\ell \le n^*+1$ a set
$B_\ell \subseteq \alpha$ decreasing with $\ell$ such that:
\mr
\item "{$(\alpha)$}"  $\{\lambda_{i,k} \in {\frak d}:i \in B_\ell
\text{ and } k \ge \ell\} \in I^+_1$
\sn
\item "{$(\beta)$}"  for each $k < \ell$ the set $\{\lambda_{i,k}:i \in
B_\ell\}$ belongs to $I_1$.
\ermn
For $\ell = 0$, the set $B_0 = \{i < \alpha:|{\frak b}_i| = n^*\}$ is
O.K.: in clause $(\alpha)$ we ask $\dbcu_{i < \alpha} {\frak b}_i \cap
{\frak d} \in I^+_1$, by which we mean ${\frak d} \in I^+_1$ which is 
assumed and Clause $(\beta)$ is empty (no $k < \ell$!) lastly by the choice
of $n^*$ we are done. \nl
For $\ell +1$, if $\ell,B_\ell$ are not as required, then there is $B'
\subseteq B_\ell$ such that

$$
\dbcu_{i \in B'} {\frak b}_i \in I_1 \text{ and } \{\lambda_{i,\ell}:i
\in B'\} \in I_1 \text{ have different truth values}.
$$
\mn
By obvious monotonicity this means $\dbcu_{i \in B'} {\frak b}_i \notin
I_1,\{\lambda_{i,\ell}:i \in B'\} \in I_1$ so let $B_{\ell + 1} = B'$. \nl
If $B_{n^*+1}$ is well defined we have by clause $(\alpha)$ that
$\{\lambda_{i,k}:i \in B_{n^*+1} \text{ and } k \ge n^* +1\} \in I^+_1$ but
as $B_{n^*+1} \subseteq B_0$ this set is empty, easy contradiction. \nl
2) Same proof except that, for defining ${\frak b}_i$, instead of quoting
\scite{2.3B} we use \nl
\cite[6.7A(h)$^*$]{Sh:430}.  We could have used it in the
proof of part (1) here. \nl
3) We apply \cite[6.7A]{Sh:430} to ${\frak a} =: \{\tau^i_j:i <
\alpha,j < \beta\} \cup \{\chi_i:i < \alpha\}$ and \wilog \, $\langle \chi_i:
i < \alpha \rangle,I,J,\langle \sigma_j:j < \beta \rangle$ and $\left <
\langle \tau^i_j:j < \beta \rangle:i < \alpha \right>$ belong to $N_0$.
Let ${\frak a}^* \in J_{\le \chi}[{\frak a}]$ be such that $J_{\le \chi}
[{\frak a}] = J_{< \chi}[{\frak a}] + {\frak a}^*$ and let ${\frak e}_j =
\{\tau^i_j:i < \alpha\} \cap {\frak a}^*$ but if possible ${\frak a}^* =
{\frak a}$.  Again by \cite[6.7A,(h)$^+$]{Sh:430} for each $j$ there is
${\frak c}_j \subseteq \text{ pcf}_{\sigma_j\text{-complete}}({\frak e}_j)$
such that ${\frak e}_j \subseteq \dbcu_{\theta \in {\frak c}_j}
{\frak b}^{\beta+1}_\theta[\bar{\frak a}]$.  Let ${\frak c}_j = \{\lambda^j
_\varepsilon:\varepsilon < \varepsilon_j\}$ with no repetitions and let
${\frak b}^j_\varepsilon = {\frak b}^{\beta +1}_{\lambda^j_\varepsilon}
[\bar{\frak a}] \cap {\frak e}_j$.

Now clause (a) holds by the choices of ${\frak c}_j$ and ${\frak b}^j
_\varepsilon$.  As for clause (b) and max pcf$({\frak b}^j_\varepsilon) = 
\lambda^j_\varepsilon$ by \scite{2.3A}, i.e. clause (j) of \relax
\cite[6.7A]{Sh:430} and clearly $\lambda^j_\varepsilon \in 
\text{ pcf}_{\sigma_j\text{-complete}}({\frak e}_j)$ but 
$\lambda^j_\varepsilon \notin \text{ pcf}({\frak e}_j \backslash 
{\frak b}^j_\varepsilon)$
by clause (e) of \cite[6.7A]{Sh:430} so necessarily $\lambda^j_\varepsilon \in
\text{ pcf}_{\sigma_j\text{-complete}}({\frak b}^j_\varepsilon)$. \nl
Let $J^*$ be the ideal with domain $\dbcu_{j < \beta} \{j\} \times
\varepsilon_j$ defined by $J^* = \{A \subseteq \text{ Dom}(J^*):
\text{max pcf}\{\lambda^j_\varepsilon:(j,\varepsilon) \in A\} < \chi\}$.  By
transitivity of pcf, $\chi \in \text{ pcf}(\{\tau^i_j:i < \alpha,j < \beta\})$
hence by the choice of ${\frak a}^*,{\frak e}_j$ clearly $\chi = 
\text{ max pcf}(\dbcu_{j < \beta}{\frak e}_j)$.
\sn
As in the proof of clause (d) of part (1) we have
\mr
\item "{$(*)$}"  for $w \subseteq \alpha$ we have
$$
\text{max pcf}(\dbcu_{i \in w}{\frak e}_i) < \chi \Leftrightarrow
\text{ max pcf}(\dbcu_{i \in w} {\frak c}_i) < \chi.
$$
\ermn
We conclude that $\chi = \text{ max pcf}(\dbcu_{\ell < \alpha}{\frak c}_i)$
hence $J^*$ satisfies clause (d) (well maybe ${\frak c}_{i_1} \cap
{\frak c}_{i_2} \ne \emptyset$?  Remember \cite[Ch.I,\S1]{Sh:g}).
\sn
Lastly, we prove clause (e) so assume $w \in J^*$, so by the definition of
$J^*$, we have max pcf$({\frak d}) < \chi$ where ${\frak d} = \{\lambda^j
_\varepsilon:(j,\varepsilon) \in w\}$.  So by transitivity of pcf
(\cite[Ch.I,1.11]{Sh:g}) as $\chi = \text{ tcf}(\dsize \prod_{i < \alpha}
\chi_i/I)$ necessarily $B =: \{i < \alpha:\chi_i \in \text{ pcf}({\frak d})
\} \in I$.  Now for each $i \in \alpha \backslash B$ we have $\chi_i
\notin \text{ pcf}({\frak d})$ hence $\chi_i \notin \text{ pcf}({\frak d} \cap
{\frak e}_i)$ but $\dsize \prod_{j < \beta} \tau^i_j/J$ has true cofinality
$\chi_i$, so necessarily $B_i =: \{j < \beta:\tau^i_j \in {\frak d} \cap
{\frak e}_i\} \in J$.  Checking the meaning you get clause (e).
\hfill$\square_{\scite{2.3}}$
\bigskip

\demo{\stag{2.4} Observation}  If $\kappa > \aleph_0,\lambda \in 
\text{ pcf}_{\kappa-\text{complete}}({\frak a})$ \ub{then} 
for some $\theta,\kappa \le
\theta = \text{ cf}(\theta) \le |{\frak a}|$, and $\langle \chi_i:i < 
\theta \rangle$ we have: $\chi_i$ regular, $\chi_i \in \lambda \cap 
\text{ pcf}({\frak a})$ and for some $\theta$-complete ideal $I \supseteq
J^{\text{\rm bd\/}}_\theta$ we have 
$\lambda = \text{ tcf}(\dsize \prod_{i < \theta} \chi_i/I)$.
\enddemo
\bigskip

\demo{Proof}  Without loss of generality
$\lambda = \text{ max pcf}({\frak a})$, otherwise replace it by
${\frak b}_\lambda[{\frak a}]$; let $J$ be the $\kappa$-complete filter
on ${\frak a}$ which $J_{<\lambda}[{\frak a}]$ generates.
Let $\theta$ be minimal such that $J$ is not
$\theta^+$-complete so necessarily $\kappa \le \theta = \text{ cf}(\theta)
\le |{\frak a}|$; as we can replace ${\frak a}$ by any ${\frak a}' 
\subseteq {\frak a},{\frak a}' \notin J_{<\lambda}[{\frak a}]$ without
loss of generality ${\frak a}$ is the union of $\theta$ members of $J$, 
so for some  ${\frak a}_i \in J$ (for $i < \theta$) we have  
${\frak a} = \dsize \bigcup_{i<\theta}
{\frak a}_i$, as $J$ is $\theta$-complete without loss of generality
${\frak a}_i \in J_{< \lambda}[{\frak a}]$.  By \scite{2.3}(1), we 
have $\langle {\frak b}_i:i < \theta \rangle,{\frak b}$ and $I$ as there.  
As $J$ is $\theta$-complete
$\{\dbcu_{i \in w} {\frak b}_i:|w| < \theta\} \subseteq I$, so by applying
clause (f), we can finish.  \hfill$\square_{\scite{2.4}}$
\enddemo
\bigskip

\demo{Proof of 2.1}  We shall prove $\otimes^1_{\lambda,\mu}$
by induction on $\lambda$.  Arriving to  $\lambda$, assume it is a
counterexample so necessarily  $\lambda > \mu$, cf$(\lambda) = 
\text{ cf}(\mu)$.  For each $\kappa < \mu$ there is ${\frak a} \subseteq
(\mu,\lambda)$ such that $|{\frak a}| < \mu$ and 
pcf$_{\kappa\text{-complete}}({\frak a}) \nsubseteq \lambda$, so by 
\cite[ChIX,4.1]{Sh:g} 
without loss of generality for some $\kappa$-complete 
ideal $J$ on ${\frak a},\lambda^+ = \text{ tcf}(\Pi {\frak a}/J)$.  So 
(by \scite{2.4}) the following subset of (cf$(\mu),\mu) \cap \text{ Reg}$ 
is unbounded in $\mu$ (by \scite{2.4}):
$$
\align
R =: \biggl\{ \theta:&\text{cf}(\mu) < \theta = \text{ cf}(\theta) 
< \mu \text{ and there is }
\langle \chi_{\theta,\zeta}:\zeta < \theta \rangle, \\
  & \text{ a sequence of regular cardinals }
\in(\mu,\lambda) \\
  & \text{ and a } \theta-\text{complete ideal } I_\theta \text{ on } \theta
\text{ extending }  J^{\text{\rm bd\/}}_\theta \text{ such that } \\
  &\underset {\zeta < \theta}\to \prod \chi_{\theta,\zeta}/I_\theta
\text{ has true cofinality } \lambda^+ \biggr\}.
\endalign
$$
\medskip
\noindent
Let $\theta,\theta_1,R^*$ be witnesses for $\otimes^0_\mu$ (i.e.
$\otimes^0_{R^*,\theta,\theta_1}$ holds),
without loss of generality otp$(R^*) = \text{ cf}(\mu)$ and remember 
cf$(\mu) < \theta_1,\theta^+ < \text{ Min}(R^*),\theta \in R$.
Let $\alpha^\ast = \theta$, we now define by induction on $n$ the 
following:  $J_n,w^n,\langle w^n_i:i < \theta \rangle,\langle \lambda_x:
x \in w^n \rangle,h^n$ as in observation \scite{2.2} such that  
$\{ x \in w^n:\lambda_x \le \mu^+ \} \in J_n$ and  $h^n_i(y) = x 
\and \lambda_x > \mu^+ \Rightarrow \lambda_y < \lambda_x$, so we 
shall get a contradiction (as we can first 
throw away the set $\{x \in  w^n:\lambda_x \le \mu^+\})$.  We also demand 
$\dsize \prod_{x \in w^n} \lambda_x/J_n$ is $\lambda^+$-directed and 
$[x \in w^n \Rightarrow \lambda_x < \lambda]$.
We let $w^0_i = \{i\},\lambda_i = \chi_{\theta,i}$, and $J_0 = I_\theta$.  
Suppose all have been defined for $n$.  Now by the induction hypothesis on 
$\lambda$ (as $\mu = \text{ sup}(R^*)$) for every $x \in w_n$,
if $\lambda_x > \mu^+$ then for some $\sigma = \sigma[\lambda_x] \in R^*$ 
we have

$$
{\frak a} \subseteq (\mu,\lambda_x) \and |{\frak a}| < \mu \Rightarrow
\text{ pcf}_{\sigma\text{-complete}}({\frak a}) \subseteq \lambda_x.
$$
\medskip
\noindent
Remember $J_n$ is $|R^*|^+$-complete (as $\theta > \text{ cf}(\mu)$), so 
it is enough to deal separately with each  $u^{n,\sigma } = u(n,\sigma) =: 
\{x \in w^n:\sigma [\lambda_x] = \sigma \text{ and } \lambda_x > \mu^+\}$ 
where $\sigma \in R^*$.
If $u^{n,\sigma} \in J_n$ we have nothing to do.  Otherwise choose  
$\kappa_\sigma \in R^\ast$, $\kappa_\sigma > \sigma,\theta$ and  
$I_{\kappa_\sigma}$, $\langle \chi_{\kappa_\sigma,\zeta}:\zeta < \kappa_\sigma
\rangle$ witnessing $\kappa_\sigma \in R$.  By \cite[Ch.IX,4.1]{Sh:g} applied 
to $\chi_{\kappa_\sigma,\zeta} < \lambda^+ = \text{ tcf} \dsize \prod
_{x \in u(n,\sigma)} \lambda_x/J_n$, for each $\zeta < \kappa_\sigma$ we can
find a sequence $\langle \tau^{n,\sigma,\zeta }_x:x \in u^{n,\sigma}\rangle$,
$\tau^{n,\sigma,\zeta}_x$ regular $< \lambda_x$ but $\ge \mu^+$ and 
$\dsize \prod_{x \in u(n,\sigma)} \tau^{n,\sigma,\zeta}_x/J_n$ has 
true cofinality $\chi_{\kappa_\sigma,\zeta}$.

Now apply \scite{2.3}(3) with $\alpha,\beta,I,J,\chi,\langle \chi_i:i <
\alpha \rangle,\langle \tau^i_j:j < \beta \rangle,\langle \sigma_j:j < \beta
\rangle$ there standing for $\kappa_\sigma,u(n,\sigma),I_{\kappa_\sigma},
J_n \restriction u(n,\sigma),\lambda^+,\langle \chi_{\kappa_\sigma,\zeta}:
\zeta < \kappa_\sigma \rangle,\langle \tau^{n,\sigma,\zeta}_x:x \in
u(n,\sigma) \rangle,\langle \sigma:x \in u(n,\sigma) \rangle$. This gives us
objects $\langle {\frak b}^{n,\sigma,\varepsilon}_x:x \in u(n,\sigma),
\varepsilon < \varepsilon_x \rangle$ and $J^{n,\sigma}$ as there.  We could
have changed some values of $\tau^{n,\sigma,\zeta}_x$ to $\mu^+$ to guarantee
that $\lambda^+ = \text{ max pcf}\{\tau^{n,\sigma,\zeta}_x:x \in u(n,\sigma),
\zeta < \kappa_\sigma\}$, so \wilog \, $\{\tau^{n,\sigma,\zeta}_x:\zeta <
\kappa_\sigma\} = \dbcu_{\varepsilon < \varepsilon_x} {\frak b}^{n,\sigma,
\varepsilon}_x$.  By \scite{2.3}(3), we have
\mr
\item "{$(*)_1$}"  if $w \subseteq \text{ Dom}(J^{n,\sigma})$ and \nl
$\{\zeta < \kappa_\sigma:\{x \in u(n,\sigma):(\exists \varepsilon <
\varepsilon_x)[\tau^{n,\sigma,\zeta}_x \in {\frak b}^{n,\sigma,\varepsilon}_x
\and (x,\varepsilon) \in w]\} \notin J_n\} \notin I_{\kappa_\sigma}$, 
then $w \notin J^{n,\sigma}$.
\ermn
Let $I^{n,\sigma}$ be the ideal on Dom$(J^{n,\sigma})$ defined by

$$
\align
w \in I^{n,\sigma} \Leftrightarrow \bigl\{ \zeta < \kappa_\sigma:\{x \in
u(n,\sigma):&(\exists \varepsilon < \varepsilon_x)
[\tau^{n,\sigma,\varepsilon}_x \in {\frak b}^{n,\sigma,\varepsilon}_x \and \\
  &(x,\varepsilon) \in w]\} \notin J_n \bigr\} \in I_{\kappa_\sigma}.
\endalign
$$
\mn
Now $(*)_1$ tells us that $J^{n,\sigma} \subseteq I^{n,\sigma}$.  Note that
since $I_{\kappa_\sigma}$ and $J_n$ are $\theta$-complete proper ideals --
we assumed $u(n,\sigma) \notin J_n$ -- we have that $I^{n,\sigma}$ is a
$\theta$-complete proper ideal on Dom$(J^{n,\sigma})$.  This means that if
we want to verify that a set is not in the $\theta$-complete ideal generated
by $J^{n,\sigma}$, it suffices to see it is not in $I^{n,\sigma}$.
\mn
By \scite{2.3}(3), clause (b) we have $\lambda^{n,\sigma,\varepsilon}_x =:
\text{ max pcf}({\frak b}^{n,\sigma,\varepsilon}_x)$ is in
pcf$_{\sigma\text{-complete}}({\frak b}^{n,\sigma,\varepsilon}_x)$.  Since
${\frak b}^{n,\sigma,\varepsilon}_x \subseteq \lambda_x$, our choice of
$\sigma$ guarantees
\mr
\item "{$(*)_2$}"  $\lambda^{n,\sigma,\varepsilon}_x = \text{ max pcf}
({\frak b}^{n,\sigma,\varepsilon}_x) < \lambda_x$.
\ermn
For $\zeta < \kappa_\sigma$, let $f^{n,\sigma}_\zeta:u(n,\sigma) \rightarrow 
\sigma$ be defined by $f_\zeta(x) = \text{ Min} \{ \epsilon < \epsilon_x:
\tau^{n,\sigma,\zeta}_x \in {\frak b}^{n,\sigma,\epsilon}_x\}$.  Now we 
can apply the choice of $\theta_1,\theta$ (i.e. for them 
$\otimes^0_{R^\ast,\theta,\theta_1}$ holds) only instead of ``$J$ a 
$\theta$-complete ideal on $\theta"$  we have 
here  $``J_n$ is a $\theta$-complete ideal on a set of cardinality $\theta$
and actually use $J_n \restriction u^{n,\sigma}"$.
So we get $A^{n,\sigma} \in I^+_{\kappa_\sigma}$ and  
$B^{n,\sigma}_\zeta = u(n,\sigma) \text{ mod } J_n$ for 
$\zeta \in A^{n,\sigma}$ such that:

$$ 
x \in  u^{n,\sigma} \Rightarrow \theta_1 > |\{ f^{n,\sigma}_\zeta(x):\zeta
\in A^{n,\sigma},x \in B^{n,\sigma}_\zeta \}|. \tag"{$(*)_3$}"
$$
\medskip
\noindent
Let us define

$$
w^{n+1}_{i,\sigma}  = \{(x,\sigma,\epsilon):(\exists \zeta \in A^{n,\sigma})
[x \in B^{n,\sigma}_\zeta \and \epsilon = f^{n,\sigma}_\zeta(x) 
\and x \in w^n_i]\}
$$

$$
h^n_{i,\sigma}:w^{n+1}_{i,\sigma} \rightarrow w^n_i \text{ is }
h^n_{i,\sigma} ((x,\sigma ,\epsilon )) = x
$$

$$
x \in u^{n,\sigma} \Rightarrow \lambda_{(x,\sigma,\epsilon)} = 
\lambda^{n,\sigma,\varepsilon}_x.
$$
\mn
Recall we are assuming $u^{n,\sigma} \in J^+_n$, if $i \in u^{n,\sigma} \in
J_n$ we let $w^{n+1}_{i,\sigma} = \emptyset$.
Now we switch ``integrating" on all $\sigma \in R^*$:

$$
w^{n+1}_i = \dbcu_{\sigma \in R^*} w^{n+1}_{i,\sigma}
$$
\mn
We let

$$
w^{n+1} = \dsize \bigcup_{\sigma \in R^*} \,\,
\dbcu_{i < \theta} w^{n+1}_{i,\sigma},h^n = 
\dsize \bigcup_{\sigma \in R^*} \,\, \dbcu_{i < \theta} h^n_{i,\sigma}.
$$

$$
\align
J_{n+1} = \biggl\{ u \subseteq w^{n+1}:&\text{ for some } i < \theta
\text{ and } u_j \subseteq u \text{ for } j < i \text{ we have } \\
 &\,u = \dsize \bigcup_{i<j}u_j \text{ and for each } j<i \text{ we have } \\
 &\,\lambda^+ > \text{ max pcf} \{\lambda_{(x,\sigma,\epsilon)}:
(x,\sigma,\epsilon) \in u_j\} \biggr\}.
\endalign
$$
\mn
Most of the verification that $w^{n+1},h^n$ and $J_{n+1}$ are as required
is routine; we concentrate on a few important points
\mr
\item "{$\boxtimes_0$}"  $|w^{n+1}_i| < \theta_1$ \nl
[Why?  By $(*)_3$, as cf$(\mu) < \theta_1 < \theta$.]
\sn
\item "{$\boxtimes_1$}"  if $x \in w^n,\lambda_x > \mu^+$ and $h^n(y) = x$,
then $\lambda_y < \lambda_x$ \nl
[Why?  Choose $\sigma$ such that $x \in u(n,\sigma)$.  If $u(n,\sigma) \in
J_n$ then $\lambda_y = \mu^+ < \lambda_x$.  If $u(n,\sigma) \notin J_n$ then
we are done by $(*)_2$.]
\sn
\item "{$\boxtimes_2$}"  $w^{n+1} \notin J_{n+1}$ \nl
[Why?  Choose $\sigma \in R^*$ with $u(n,\sigma) \notin J_n$, and let
$v(n,\sigma) = \{(x,\varepsilon):(x,\sigma,\varepsilon) \in w^{n+1}_\sigma\}$.
\nl
For $\zeta \in A^{n,\sigma}$,
$$
B^{n,\sigma}_\zeta \subseteq \{x \in u(n,\sigma):(\exists \varepsilon <
\varepsilon_x)[\tau^{n,\sigma,\zeta}_x \in {\frak b}^{n,\sigma,\varepsilon}_x
\wedge (x,\varepsilon) \in v(n,\sigma)]\},
$$
and so $v(n,\sigma) \notin I^{n,\sigma}$.  Thus $v(n,\sigma)$ is not in the
$\theta$-complete ideal generated by $J^{n,\sigma}$, and the definitions of
$J^{n,\sigma}$ and $J_{n+1}$ imply $w^{n+1}_\sigma \notin J_{n+1}$.]
\sn
\item "{$\boxtimes_3$}"  For every $A \in J_{n+1},\{x \in w^n:(\forall y \in
w^{n+1})[h^n(y) =x \Rightarrow y \in A]\} \in J_n$. \nl
[Why?  Suppose $B \in J^+_n$, and choose $\sigma \in R^*$ such that $B \cap
u(n,\sigma) \in J^+_n$.  Let $A = \{(x,\sigma,\varepsilon) \in
w^{n+1}:x \in B\}$, and let $A' = \{(x,\varepsilon):(x,\sigma,\varepsilon)
\in A\}$.  For $\zeta \in A^{n,\sigma}$
$$
B \cap B^{n,\sigma}_\zeta \subseteq \{x \in u(n,\sigma):(\exists \varepsilon
< \varepsilon_x)[\tau^{n,\sigma,\zeta}_x \in {\frak b}^{n,\sigma,\varepsilon}
_x \wedge (x,\varepsilon) \in A']\}, 
$$
and since $B \cap B^{n,\sigma}_\zeta \in J^+_n$, we know $A' \notin
I^{n,\sigma}$ hence $A \notin J_{n+1}$.]
\ermn
Thus we have carried out the induction and hence get by \scite{2.2} 
the contradiction and finish the proof.  \hfill$\square_{\scite{2.1}}$
\enddemo
\bigskip

\remark{\stag{2.4A} Remark}  1) We can be more specific phrasing \scite{2.1}:
let $R^\ast \subseteq \mu$  be unbounded, ${\bar \Gamma} =
\langle \Gamma_\sigma :\sigma \in R^\ast \rangle$,
$\Gamma_\sigma$ a set of ideals on $\sigma$; the desired
conclusion is: for every  $\lambda  > \mu$ for some $\sigma^\ast < \mu$  
we have: if  $\sigma \in R^* \backslash \sigma^\ast$,  $\lambda_i \in
(\mu,\lambda) \cap \text{ Reg}$ for $i < \sigma$, $J$, $J \in \Gamma_\sigma$ 
then $\text{pcf}_{\Gamma_\sigma}\left(\dsize \prod_{i<\sigma} \lambda_i,
\le_J \right) \subseteq \lambda$.  (Reg is the class of regular cardinals).\nl
2)  You can read the proofs for the case $\mu$ strong limit singular 
and get an alternative to the proof in {\S1}.
\endremark
\bigskip

\proclaim{\stag{2.5} Claim}  Assume $\lambda^* > \mu > \aleph_1,\mu$ an 
uncountable limit cardinal and we have: 
\roster
\item "{$\otimes^{1.5}_{\lambda^\ast,\mu}$}"  for every
$\lambda \in (\mu,\lambda^\ast]$, we have $\otimes^1_{\lambda,\mu}$
 (from the conclusion of 2.1).
\endroster
\ub{Then}
\roster
\item "{$\otimes^2_{\lambda^\ast,\mu}$}"  $(\alpha)$\,\,${\frak a}
 \subseteq (\mu,\lambda^\ast), {\frak a} \subseteq \text{{\rm Reg\/}},
|{\frak a}| < \mu \Rightarrow |\lambda^* \cap \text{{\rm pcf\/}}
({\frak a})| \le \mu$
\sn
\item "{${}$}"  $(\beta)$\,\, if $\mu$ is regular then (for 
${\frak a} \subseteq \text{{\rm Reg\/}}$):
\endroster
$$
{\frak a} \subseteq (\mu,\lambda^\ast) \and |{\frak a}| < \mu
 \Rightarrow |\lambda^\ast \cap \text{{\rm pcf\/}}({\frak a})| < \mu.
$$
\endproclaim
\bigskip

\demo{Proof}  Let
$\mu(*) = 
\left\{
\alignedat2
&\mu^+ &&\text{ if $\mu$ is singular } \\
&\mu &&\text{ if $\mu$ is regular }.\endalignedat
\right.$
\medskip

So assume ${\frak a} \subseteq (\mu,\lambda^*) \cap \text{{\rm Reg\/}},
|{\frak a}| < \mu$, and $\lambda^\ast \cap 
\text{ pcf}({\frak a})$ has cardinality $\ge \mu(*)$.  
Let $\lambda_0 = \text{ Min}({\frak a})$ and $\langle \lambda_{i+1}:
i \le \mu(*) \rangle$ list the first $(\mu(*) + 1)$-members of 
$\lambda^* \cap (\text{pcf}({\frak a})) \backslash \{ \lambda_0 \}$  
(remember $\text{ pcf}({\frak a})$ has a last member), and for limit
$\delta \le \mu(*)$, let  $\lambda_\delta = \dsize \bigcup_{i<\delta}
\lambda_i$ so $\lambda_{\mu(*)} \le \lambda^\ast$.  Now by an assumption
for some  $\kappa < \mu$, $\otimes^1_{\lambda_{\mu(*)},\mu,\kappa}$ (from
2.1), without loss of generality  
$\kappa$ is regular.  Now choose by induction on  $\zeta < \mu$,
$i(\zeta)$ such that $i(\zeta) < \mu(*)$ is a successor ordinal,
$i(\zeta) > \dsize \bigcup_{\xi < \zeta} i(\xi)$, and
$\lambda_{i(\zeta)} > \sup{\text{pcf}_{\kappa \text{-complete}}}
( \{ \lambda_{i(\xi)}:\xi < \zeta \})$. \newline
\noindent
Why is this possible?  We know  
${\text{pcf}_{\kappa \text{-complete}}} ( \{ \lambda_{i(\xi)}:\xi < \zeta \})$
cannot have
a member  $\ge \lambda_{\mu(*)}$,  (hence  $> \lambda_{\mu(*)}$ being 
regular), by the choice of  $\kappa$.  Also \newline
${\text{pcf}_{\kappa \text{-complete}}}
( \{ \lambda_{i(\xi)}:\xi < \zeta \})$ cannot be unbounded in
$\lambda_{\mu(*)}$ (because cf$(\lambda_{\mu(*)}) = \mu(*) \ge \kappa$
(remember $\mu(*)$ is regular) as then it will have a member 
$> \lambda_{\mu(*)}$, see \cite[Ch.I,1.11]{Sh:g}).  
So it is bounded below $\lambda_{\mu(*)}$ hence $i(\zeta)$ exists. 

Now we get contradiction to \cite[3.5]{Sh:410}, version (b) of (iv) 
there \newline
(use e.g. $\langle \lambda_{i(\zeta)}:\zeta < (\kappa + |{\frak a}|)^{+4}
\rangle$). (Alternatively to \cite[6.7F(5)]{Sh:430}).
\hfill$\square_{\scite{2.5}}$
\enddemo
\bigskip

\proclaim{\stag{2.6} Theorem}  Let  $\mu$  be a limit uncountable singular
cardinal, $\mu < \lambda$ and \newline
$[|{\frak a}| < \mu \and {\frak a} \subseteq \text{{\rm Reg\/}} \cap
(\mu,\lambda) \Rightarrow  |\lambda \cap \text{{\rm pcf\/}}
({\frak a})| < \mu]$
\newline
or at least:

$$
\text{ for every large enough }  \sigma \in \text{{\rm Reg\/}} \cap \mu,
\text{ we have:} \tag"{$\oplus_{\mu,\lambda}$}"
$$

$$
\text{\underbar{if} } \, {\frak a} \subseteq \text{{\rm Reg\/}} \cap 
(\mu,\lambda),|{\frak a}| < \mu \,\, \text{\underbar{then}} \,\,
|\lambda \cap {\text{{\rm pcf\/}}_{\sigma \text{-complete}}}({\frak a})|< \mu
\tag"{$\oplus^\sigma_{\mu,\lambda}$}"
$$  
\mn
\underbar{Then} for every large enough  $\kappa < \mu$ we have
$\otimes^1_{\mu,\kappa}$ of \scite{2.1}, hence 
$\text{{\rm cov\/}}(\lambda,\mu,\mu,\kappa) = \lambda$.
\endproclaim
\bigskip

\remark{Remark}  This proof relies on \cite[\S5]{Sh:420}.
\endremark
\bigskip

\demo{Proof}  Without loss of generality cf$(\mu) = \aleph_0$ (e.g.
force by $\text{Levy}(\aleph_0$,cf$(\mu))$ as nothing relevant changes, 
or argue as in \scite{1.2A}, however, we can just repeat the proof). 

Assume this fails.  Without loss of generality $\lambda$ is minimal, so  
cf$(\lambda) = \aleph_0$.  
Failure means (by \scite{2.4}) that $\mu = \sup(R)$ when

$$
\align
R = \biggl\{ \theta:&\theta \in \mu \cap \text{ Reg and for some }
\chi_\zeta \in \text{ Reg } \cap (\mu,\lambda) \text{ for } \zeta < \theta, \\
  &\text{ and } \theta{\text{-complete ideal}} \,\, I \text{ on } \theta \\
  &\text{ we have }
\lambda^+ = \text{ tcf}(\dsize \prod_{\zeta < \theta} \chi_\zeta/I)\biggr\}.
\endalign
$$
\mn
For simplicity assume that for $\chi < \mu$ and $A \subseteq (2^\chi)^+$,
in $K[A]$ there are Ramsey cardinals $> \chi$.  This makes a minor restriction
say for one $\lambda$ we may get $\le \lambda^+$ instead of $< \lambda^+$
(which is equivalent to $< \lambda$). \nl
So by \cite[\S5]{Sh:420}, for some uncountable regular $\sigma < \kappa$
from $R \backslash \text{ cf}(\mu)^+,\oplus^\sigma_{\mu,\lambda}$ from the
assumption of the theorem holds and for some family $E$ of ideals on 
$\kappa$ normal by a function $\imath:\kappa \rightarrow \sigma$ and
$J \in E$ and $\lambda_i = \text{ cf}(\lambda_i) \in (\mu,\lambda),
\lambda^+ = \text{ tcf} 
\left( \dsize \prod_{i < \kappa} \lambda_i/J \right)$ and
$\langle \lambda_i:i < \kappa \rangle$,
$J$  minimal in a suitable sense, that is $\alpha(*) = \text{ rk}^3_J(\langle 
\lambda_i:i < \kappa \rangle,E$) is minimal
so without loss of generality rk$^3_J( \langle \lambda_i:i < \kappa \rangle,
E) = \text{ rk}^2_J (\langle \lambda_i:i < \kappa \rangle,E$).  Hence we 
do not have  $A \subseteq \kappa$, $\kappa \backslash A \notin J$ and
$\lambda^\prime_i \in (\mu,\lambda) \cap$ Reg such that
$\langle \lambda^\prime_i:i < \kappa \rangle <_{J+A} \langle \lambda_i:
i < \kappa \rangle$ and $\lambda^+ = \text{ tcf}(\dsize \prod_{i < \kappa}
\lambda^\prime_i/J)$.  As cf$(\mu) = \aleph_0$, we can find  
$\langle \theta_n:n < \omega \rangle,\kappa < \theta_n \in R \cap  
\mu$ and $\mu = \dsize \bigcup_{n < \omega} \theta_n$.  As
$\lambda$ is minimal there is a partition $\langle u(n):n < \omega \rangle$
of $\kappa$, such that:

$$
i \in u(n), n < \omega, |{\frak a}| < \mu, {\frak a} \subseteq \text{ Reg } 
\cap (\mu,\lambda_i) \Rightarrow 
{\text{pcf}_{\theta_n \text{-complete}}}({\frak a})
\subseteq \lambda_i. \tag"{$(*)$}"
$$
\mn
So for some $n$ we have $u(n) \in J^+$.
Without loss of generality $(\forall i < \kappa)(\lambda_i > \mu^+)$ and
 (as $\sigma > \aleph_0)$ for some  $n = n(\ast)$ we have $u(n) = \kappa$
(i.e. the minimality of $\alpha(*)$
is preserved).  Choose $\theta \in R \cap \mu$ large enough such that 
$(\forall{\frak a}) \biggl[{\frak a} \subseteq \text{ Reg } \cap 
(\mu,\lambda) \and |{\frak a}| \le \theta_{n(*)} + \kappa
\Rightarrow |\lambda \cap {\text{ pcf}_{\sigma \text{-complete}}}({\frak a})|
< \theta \biggr]$.  (Why is this possible?  As $\oplus^\sigma_{\mu,\lambda}$
which holds by the choice of $\sigma$).
As $\theta \in R \cap \mu$ we can choose a sequence
$\langle \chi_\zeta:\zeta < \theta \rangle$ and 
$I \supseteq J^{\text{\rm bd\/}}_\theta$ a 
$\theta$-complete ideal on $\theta$ such that
$\chi_\zeta \in (\mu,\lambda)$ and tcf$\left( \dsize \prod_
{\zeta < \theta}\chi_\zeta/I \right) = \lambda^+$.  By \cite[Ch.IX,4.1]{Sh:g}
we can find $\tau^\zeta_i = \text{ cf}(\tau^\zeta_i) \in (\mu,\lambda_i),
\tau^\zeta_i < \lambda_i$ such that $\chi_\zeta = \text{ tcf}
\left( \dsize \prod_{i < \kappa} \tau^\zeta_i/J \right)$.

Now  ${\frak a} =: \lambda \, \cap {\text{ pcf}_{\sigma \text{-complete}}}
\{ \tau^\zeta_i:i < \kappa,\zeta < \theta \}$ has cardinality $< \mu$ 
(by the choice of $\sigma$) and has a smooth closed representation 
$\langle {\frak b}_\Upsilon({\frak a}):\Upsilon \in {\frak a} \rangle$ 
(see \cite[6.7]{Sh:430}).  For $i < \kappa$ there is 
${\frak c}_i \subseteq \text{ pcf}_{\theta_{n(*)}\text{-complete}}
\{ \tau^\zeta_i:\zeta < \theta\}$ such that 
$|{\frak c}_i| < \theta_{n(*)}$ and 
$\dsize \bigwedge_{\zeta < \theta} \tau^\zeta_i \in \cup
\{ {\frak b}_\Upsilon({\frak a}):\Upsilon \in {\frak c}_i \}$
(by the choice of $n(*)$ and by \cite[6.7]{Sh:430}, note that $\sigma <
\kappa < \theta_{n(*)}$ by their choices hence
pcf$_{\theta_{n(*)}\text{-complete}}\{\tau^\zeta_i:\zeta < \theta\}
\subseteq {\frak a}$ hence all is O.K.). 
Also ${\frak c}_i \subseteq \lambda_i$ because 
we are assuming $u_{n(*)} = \kappa$.

Let
$$
\align
{\frak d} =: \biggl\{ \text{tcf}(\dsize \prod_{i \in A} 
\tau_i/(J + A)):&\langle
  \tau_i:i < \kappa \rangle \in \dsize \prod_{i < \kappa}{\frak c}_i 
 \text{ and } A \in J^+ \text{ and } \\
  &\text{tcf}(\dsize \prod_{i \in A} \tau_i /(J + A))
 \text{ is well defined } \biggr\}
\endalign
$$
\medskip
\noindent
Let ${\frak c} = \dsize \bigcup_{i < \kappa}{\frak c}_i$.  So $|{\frak c}|
\le \kappa + \theta_{n(*)}$ hence $\lambda \cap {\text{ pcf}_{\sigma \text
{-complete}}}({\frak c})$ has cardinality $< \theta$, and ${\frak d} \subseteq
\lambda$ by the choice of $\alpha(*)$ and 
${\frak d} \subseteq {\text{ pcf}_{\sigma \text{-complete}}}
({\frak c})$ hence $|{\frak d}| < \theta$ (by the choice of $\theta$).
 
Now if $\psi \in \lambda^+ \cap \text{ pcf}({\frak c})$ then 

$$
B_\psi = \{ \zeta < \theta: \{i < \kappa:\tau^\zeta_i \in {\frak b}_\psi
[{\frak a}]\} \notin J \} \in I.
$$
\mn
[Why?  Otherwise $\zeta \in B_\psi \Rightarrow \chi_\zeta \in \text{ pcf}
({\frak b}_\psi[{\frak a}])$ hence $\text{pcf}({\frak b}_\psi[{\frak a}])$
includes \nl
pcf$\{ \chi_\zeta:\zeta \in B_\psi \}$, but as $B_\psi \notin I$ the
cardinal $\lambda^+$ belongs to the latter; but max pcf$({\frak b}_\psi
[{\frak a}]) = \psi < \lambda$ contradiction]. \newline
\noindent
But we know that $|{\frak d}| < \theta$,  
and $I$ is $\theta$-complete and ${\frak d} \subseteq \text{ pcf}({\frak c})$,
so

$$
\align
X = \biggl\{ \zeta < \theta:&\text{ for some } \psi \in  {\frak d} 
\text{ we have} \\
  &\{i < \kappa:\tau^\zeta_i \in {\frak b}_\psi[{\frak a}] \} \notin J
\biggr\} \subseteq \dbcu_{\psi \in {\frak d}} B_\psi \in I.
\endalign
$$ 
\mn
So there is some $\zeta^* \in \theta \backslash X$, 
and for $i < \kappa$ choose
$\Upsilon_i \in {\frak c}_i$ such that $\tau^{\zeta^*}_i \in 
{\frak b}_{\Upsilon_i}[{\frak a}]$ (well defined by the choice of
${\frak c}_i$). So by smoothness of the representation

$$
\psi \in {\frak d} \Rightarrow \{i < \kappa:\Upsilon_i \in {\frak b}_\psi
[{\frak a}]\} \subseteq \{i < \kappa:\tau^{\zeta^*}_i \in {\frak b}_\psi
[{\frak a}]\} \in J.
$$
\mn
Now by the pcf theorem for some $A \in J^+$ we have $\dsize \prod
_{i \in A} \Upsilon^{\zeta^*}_i/J$ has true cofinality which we call
$\Upsilon$, so necessarily $\Upsilon \in \text{ pcf}_{\sigma\text{-complete}}
(\{\Upsilon^{\zeta^*}_i:i \in A\}) \in {\frak d}$ (see the definition of
${\frak d}$) but
this contradicts the previous sentence (recall ${\frak d} \subseteq \lambda$
by the minimality of $\alpha(*)$).  \hfill$\square_{\scite{2.6}}$
\enddemo
\newpage

\head {\S3 Applications} \endhead  \resetall 
\bn
Of course
\proclaim{\stag{3.0} Claim}  If $\mu$ is as in \scite{2.1}, then the
conclusions of \scite{1.2}  and \scite{1.1} hold.
\endproclaim
\bigskip

\proclaim{\stag{3.1} Claim}  If  $\lambda \ge \beth_\omega$ then:
\roster
\item "{(a)}"  $2^\lambda  = \lambda^+ \Leftrightarrow
\diamondsuit_{\lambda^+}$
\sn
\item "{(b)}"  $\lambda = \lambda^{<\lambda}$  iff  $(D\ell)_\lambda$.
\endroster
\endproclaim
\bn
Where we remember
\definition{\stag{3.2} Definition}  1)  $(D\ell)_\lambda$ means that:
\mr
\item "{{}}"  $\lambda$ is regular uncountable and there is ${\bar \Cal P} = 
\langle {\Cal P}_\alpha:\alpha < \lambda \rangle$  such that
${\Cal P}_\alpha$ is a family of $< \lambda$  subsets of $\alpha$ satisfying:
\sn
\item "{$(*)$}"  for every  $A \subseteq \lambda$, $\{ \alpha < \lambda:
A \cap \alpha \in {\Cal P}_\alpha \}$ is a stationary subset of $\lambda$.
\ermn
2)  $(D\ell)^*_S$ \, $(S \subseteq \lambda$  stationary) means $\lambda$
regular and there is  ${\bar \Cal P}$ as above such that:
\roster
\item "{$(*)$}"  for every  $A \subseteq \lambda$ we have
$\{ \alpha \in S:A \cap x \notin {\Cal P}_\alpha \}$ is not stationary.
\endroster
\medskip
\noindent
3)  $(D\ell)^+_S$ where $S \subseteq \lambda$ is stationary, $\lambda$  
regular uncountable means that: for some  ${\bar{\Cal P}}$ as above:
\roster
\item "{$(*)$}"  for every  $A \subseteq \lambda$ for some club $C$ of  
$\lambda$  we have: \newline
$\delta \in S \cap C \Rightarrow A \cap \delta \in
{\Cal P}_\delta \and C \cap \delta \in {\Cal P}_\delta$.
\endroster 
\medskip
\noindent
4)  Let  $\lambda$  be regular uncountable,  $S \subseteq \lambda$  
stationary.  Now  $\diamondsuit_S$ means that there is \newline
$\langle A_\alpha:\alpha  
\in  S \rangle$ such that $A_\alpha \subseteq \alpha$ and for every
$A \subseteq \lambda$ the set $\{\alpha \in S:A \cap \alpha = A_\alpha \}$
is a stationary subset of  $\lambda$. \nl
5) For $\lambda$ regular uncountable and $S \subseteq \lambda$ stationary
$(D \ell)_S$ means that for some $\langle {\Cal P}_\alpha:\alpha \in S
\rangle$ as above for every $A \subseteq \lambda$ the set $\{\delta \in S:
A \cap \delta \in {\Cal P}_\delta\}$ is stationary.
\enddefinition
\bigskip

\remark{\stag{3.3} Remark}  1) If $\lambda$  is a successor cardinal,
$(D\ell)_\lambda$ is equivalent to $\diamondsuit_\lambda$
(by Kunen, so (a) is a particular case of (b) in \scite{3.1}). \newline
2)  By \cite{Sh:82}, \cite{HLSh:162}, if $(D\ell)_\lambda$ then the omitting 
types theorem for $L(Q)$ for $\lambda$-compact models in the 
$\lambda^+$-interpretation 
holds (and more).  Now  $\lambda  = \lambda^{<\lambda}$ is the standard 
assumption to the completeness theorem of  $L(Q)$  in the 
$\lambda^+$-interpretation; and is necessary and sufficient when we restrict 
ourselves to $\lambda$-compact models.  So the question arises, how strong
is this extra assumption?  If G.C.H. holds $(D\ell)_\lambda  
\Leftrightarrow \lambda = \lambda^{<\lambda}$ for every  $\lambda \ne
\aleph_1$ (by \cite{Sh:82}, continuing Gregory \cite{Gr}); and more there.  
Here we improve those theorems.  Now \scite{3.1} says that above 
$\beth_\omega$, the two conditions are equivalent. \newline
3) We may consider the function $h:\lambda \rightarrow \lambda \cap
\text{ Car}$, demanding $|{\Cal P}_\alpha| < h(\alpha)$. \nl
4) Remember that for $\lambda > \aleph_0$ regular and stationary $S = S_1
\subseteq S_2 \subseteq \lambda$ we have $(D \ell)^+_S \Rightarrow
(D \ell)^*_S \Rightarrow (D \ell)_S$ and $(D \ell)_{S_1} \Rightarrow
(D \ell)_{S_2}$ but $(D \ell)^*_S \Rightarrow (D \ell)^*_{S_1},
(D \ell)^+_{S_2} \Rightarrow (D \ell)^+_{S_1}$.
\endremark
\bigskip

\demo{\stag{3.4} Proof of 3.1}  Trivially  $(D\ell)_\lambda \Rightarrow \lambda
= \lambda^{<\lambda}$, so assume $\lambda = \lambda^{<\lambda}$,
and let \nl
$\{ A^*_i:i < \lambda \}$ list the bounded subsets of $\lambda$,
each appearing  $\lambda$ times. \newline
\noindent
For each  $\alpha < \lambda$ let

$$
R_\alpha = \{ \kappa < \beth_\omega:\text{ cov}(|\alpha|,\kappa^+,\kappa^+,
\kappa) < \lambda \text{ and } \kappa \text{ is regular}\}.
$$
\mn
We know (by \scite{1.1}) that for each $\alpha \in (\beth_\omega,\lambda)$,
$R_\alpha$ contains a co-bounded subset of $\text{ Reg } \cap \beth_\omega$,
say $\text{Reg } \cap \beth_\omega \backslash \beth_{n_\alpha}$.  So for some
$n^\ast < \omega$

$$
S^\ast = \{ \alpha < \lambda:\alpha > \beth_\omega,n_\alpha < n^\ast \}
$$
\mn
is unbounded in $\lambda$;  hence trivially  $S^* = (\beth_\omega,\lambda)$. 
So $R =: \{ \kappa < \lambda:\kappa \text{ is regular } 2^\kappa, < \lambda$ 
and for every $\alpha < \lambda$ we have cov$(|\alpha|,\kappa^+,\kappa^+,
\kappa) < \lambda \}$ contains $\text{ Reg } \cap (\beth_{n^*},\beth_\omega)$.
As $\lambda = \text{ cf}(\lambda) > \beth_\omega$, for each $\alpha <
\lambda$,  $\kappa \in R$  there is ${\Cal P}^\kappa_\alpha$, a family of 
$< \lambda$ subsets of $\alpha$ of cardinality  $\kappa$ such that if
$A \subseteq \alpha$, $|A| = \kappa$ then $A$ is included in the union
of $< \kappa$ members of ${\Cal P}^\kappa_\alpha$. \newline
\noindent
Let  ${\Cal P}^\ast_\alpha = \{ B: \text{ for some } \kappa \in R \cap
(\alpha + 1) \text{ and }  A \in {\Cal P}^\kappa_\alpha \text{ we have }
B \subseteq A\}$ so ${\Cal P}^*_\alpha$ is a family of $< \lambda$ subsets of
$\alpha$.  For each $A \subseteq \lambda$  we define
$h_A:\lambda \rightarrow \lambda$  by defining  $h_A(\alpha)$ by induction
on  $\alpha:\text{for } \alpha \text{ non-limit } h_A(\alpha)$ 
is the first ordinal  $i > \dbcu_{\beta < \alpha} h_A(\beta) +1$ such that  
$A \cap \alpha = A^*_i$ and for $\alpha$ limit $h_A(\alpha) =
\dbcu_{\beta < \alpha} h_A(\beta)$.  So $h_A(\alpha)$ is strictly increasing 
continuous, hence $h_A(\alpha) \ge \alpha$ and $h(\alpha) = \alpha
\leftrightarrow [(\alpha \text{ limit}) \and (\forall \beta < \alpha)
(h_A(\beta) < \alpha)]$.  Let

$$
{\Cal P}^0_\alpha  =: \biggl\{ \dbcu_{\beta \in B} A^*_\beta:
B \in {\Cal P}^\ast_\alpha \biggr\} 
$$

$$
{\Cal P}_\alpha =: {\Cal P}^0_\alpha \cup \biggl\{ \{ \beta < 
\alpha:h_A(\beta) = \beta \}:A \in {\Cal P}^0_\alpha \biggr\} 
$$
\medskip
\noindent
(remember $\langle A^\ast_\alpha:\alpha < \lambda \rangle$  lists the 
bounded subsets of  $\lambda$ each appearing unboundedly often). \newline
Now for any $A \subseteq \lambda$ we have
$E =: E_A =: \biggl\{ \delta < \lambda:
\delta \text{ limit and } \dsize \bigwedge_{\beta < \delta}
h_A(\beta) < \delta \biggr\}$ is a club of  $\lambda$, and
\roster
\item "{$(*)_1$}"  cf$(\delta) < \delta \in  E \and \text{ cf}(\delta) \in R
\Rightarrow A \cap \delta \in {\Cal P}^0_\delta \subseteq
{\Cal P}_\delta$
\endroster
\medskip
\noindent
[Why?  Let $\kappa =: \text{ cf}(\delta)$, and let $\langle \beta_j:j < \kappa \rangle
$ be an increasing sequence of successor ordinals with limit $\delta$,
hence $\langle h_A(\beta_j):j < \kappa
\rangle$ is (strictly) increasing with limit $\delta$; so for some $\beta
< \kappa = \text{ cf}(\delta)$ and $B_i \in 
{\Cal P}^{\text{cf}(\delta)}_\delta$ for $i < \beta$ we have 
$\{h_A(\beta_j):j < \kappa \} \subseteq
\dsize \bigcup_{i < \beta} B_i$, so for some $i$, $\{ h_A(\beta_j):
j < \kappa,h_A(\beta_j) \in B_i \}$ is unbounded in
$\delta$, and clearly $B^\prime =: \{ h_A(\beta_j):j < \kappa \} \cap
 B_i \in {\Cal P}^\ast_\delta$, hence $\cup \{ A^\ast_\gamma:\gamma \in
B^\prime \} \in {\Cal P}^0_\alpha$ is as required], and
\medskip
\roster
\item "{$(*)_2$}"  cf$(\delta) < \delta \in E \and \text{ cf}(\delta) \in R
\Rightarrow  E \cap  \delta  \in  {\Cal P}$
\endroster
\medskip
\noindent                      
[Why?  As $A \subseteq \lambda$, $\delta \in E_A \Rightarrow h_A \restriction
\delta = h_{A \cap \delta} \restriction \delta$].\hfill$\square_{\scite{3.1}}$
\enddemo
\bn
Note that we actually proved also
\proclaim{\stag{3.5} Claim}  1) Assume $\lambda = \mu^+ = 2^\mu > \chi$,
$\chi$  strong limit then for some $\chi^\ast < \chi$  we have  
$\diamondsuit^+_{ \{ \delta < \lambda:\chi^\ast < \text{cf}(\delta) 
< \chi \}}$.
\newline
\noindent
2)  Similarly for $\lambda = \lambda^{<\lambda}$ inaccessible, $\chi$  
strong limit  $< \lambda$  for some  $\chi^\ast < \chi$,
$(D\ell )^+_{ \{ \delta < \lambda:\chi^\ast < \text{ cf}(\delta) 
< \chi \} }$ holds.
\newline
\noindent
3)  If  $\lambda = \lambda^{<\lambda}$, and \newline
$S = \{ \delta < \lambda:\text{{\rm cf\/}}(\delta) < \delta,
2^{\text{{\rm cf\/}}(\delta)} < \lambda$, and 
$[\lambda > \text{{\rm cov\/}}(|\delta|,\text{{\rm cf\/}}
(\delta)^+,\text{{\rm cf\/}}(\delta)^+,\text{{\rm cf\/}}(\delta)]\}$
\underbar{then} $(D\ell)^+_S$; so if $\lambda$  is a successor 
cardinal we have $\diamondsuit^+_S$. \newline
4)  Assume \footnote{if $\lambda = \mu^+,\mu = \text{{\rm cf\/}}(\mu) > 
\theta = \text{{\rm cf\/}}(\theta) > \sigma = \text{{\rm cf\/}}(\sigma)$ 
then there are $S,\bar C$ 
as in \scite{3.5}(4)(see \cite[\S4]{Sh:351} or \cite[Ch.III,2.14]{Sh:g}).  
Of course, we get not just guessing on a stationary
set but on a positive set modulo a larger ideal.}
$\lambda = \lambda^{<\lambda} > \theta = \text{{\rm cf\/}}(\theta) > \sigma
= \text{{\rm cf\/}}(\sigma),\theta^\sigma < \lambda,S \subseteq \lambda$,
$\{ \delta \in S:\text{{\rm cf\/}}(\delta) = \theta\}$ is stationary
${\bar C} = \langle C_\alpha:\alpha \in S \rangle$, for $\alpha \in S$,
$C_\alpha$ is a closed subset of  $\alpha$,  $[\beta \in C_\alpha
\Rightarrow  \beta \in  S \and C_\beta  = \beta \cap C_\alpha]$.
Assume further that for no  $\alpha < \lambda$  is there
${\Cal P} \subseteq \{ a \subseteq \alpha:|a| = \theta \}$, such that
$[a \in {\Cal P} \and b \in {\Cal P} \and a \ne b \Rightarrow
|a \cap b| < \sigma]$, and $[{\frak a} \subseteq \lambda \cap
\text{{\rm Reg\/}}\backslash \text{{\rm Min\/}} \, 
{\frak a} > \theta \and |{\frak a}|
< \theta \Rightarrow \lambda > \sup(\lambda \cap \text{{\rm pcf\/}}
{\frak a})]$ \newline
(e.g. $\lambda$  successor). \newline
\underbar{Then} $(D\ell)_{S_\sigma}$ holds where  
$S_\sigma  = \{ \delta \in S:\text{{\rm cf\/}}(\delta) = \sigma\}$.
\endproclaim
\bigskip

\demo{Proof}  Easy, e.g.: 4) By \cite[Ch.III,\S2]{Sh:g} 
without loss of generality for every club $E$ of $\lambda$ for some  
$\delta \in E$, $C_\delta \subseteq E$.
Let $\chi = \beth_3(\lambda)^+$, let  $\langle M_i:i < \lambda \rangle$
be such that: $M_i \prec ({\Cal H}(\chi),\in,<^*_\chi),\| M_i \| < \lambda,
\lambda \in M_i,M_i \cap \lambda$ an ordinal, $\langle M_j:j \le i \rangle 
\in M_{i+1}$.  Let for 
$\delta \in S_\sigma$,  ${\Cal P}_\delta = M_{\delta +1} \cap  
{\Cal P}(\delta)$.  It is enough to show that  ${\bar{\Cal P}} = 
\langle {\Cal P}_\delta:\delta \in S_\sigma \rangle$  exemplifies 
$(D\ell)_{S_\sigma}$.  So let  $\langle x_\alpha :\alpha < \lambda \rangle  
\in  M_0$ list the bounded subsets of $\lambda$ each appearing $\lambda$  
times.  Let  $X \subseteq \lambda,E_0$ be a club of $\lambda$;  we 
define by induction on  $\alpha,h_X(\alpha) < \lambda$ as the first  
$\gamma  < \lambda$  such that  $\gamma > \dsize \bigcup_{\beta < \alpha}
h_X(\beta)$  and  $X \cap \alpha = X_\alpha$.
Let  
$\langle M^\ast_i:i < \lambda \rangle$  be chosen as above but also
$h_X \in M^\ast_0$,  $\langle M_i:i < \lambda \rangle \in M^\ast_0$,
$E_0 \in M^\ast_0$.  Let $E =: \{ \delta \in E_0:M^\ast_\delta \cap \lambda
= \delta = M_\delta \cap \lambda \}$,  clearly it is a club of $\lambda$.  
Let $\delta \in S \cap E$, cf$(\delta) = \theta$  be such that
$C_\delta \subseteq E$.  
Now we imitate the proof \cite[\S6]{Sh:410} for $h_X \restriction C_\delta$.
\nl
${{}}$  \hfill$\square_{\scite{3.5}}$
\enddemo
\bn
\proclaim{\stag{3.5A} Claim}  Above instead of demanding on $\kappa \, ``\kappa
= \text{{\rm cf\/}}(\kappa) \and 2^\kappa < \lambda \and [\alpha < \lambda
\Rightarrow \text{{\rm cov\/}}(|\alpha|,\kappa^+,\kappa^+,\kappa) < \lambda]"$
it suffices to demand ``$\kappa = \text{{\rm cf\/}}(\kappa) < \lambda$ and 
if $T$ is a tree with $\kappa$-levels and $< \lambda$ nodes then $T$ has
$< \lambda \, \kappa$-branches". \nl
See {\rm\cite[\S2\/]{Sh:589}}.
\endproclaim
\bn
\centerline {$* \qquad * \qquad *$}
\bn
\proclaim{\stag{3.6} Lemma}  1) Suppose  $c\ell$ is \ub{an operation} on 
$X$, i.e., $c\ell$ is a function from  ${\Cal P}(X)$  to  ${\Cal P}(X)$.
Assume further $\kappa \le \kappa^\ast < \mu = \mu^\kappa$ and we let

$$ 
\align
{\Cal P}^\ast  = \biggl\{
A \subseteq X:&|A| = \mu \text{ and for every } B \subseteq  A
\text{ satisfying } |B| = \kappa^\ast \text{ there is } \\
  &B^\prime \subseteq B, |B^\prime| = \kappa \text{ such that }
c\ell(B') \subseteq  A, \\
  &\text{ and } |c\ell(B')| = \mu \biggl\}.
\endalign
$$ 
\medskip
\noindent
If $\kappa^\ast < \beth_\omega(\kappa) \le \mu$ \ub{then} there is function
$h:X \rightarrow \mu$ such that: if $A \in {\Cal P}^\ast$ then
$h \restriction A$ is onto $\mu$. \newline
\noindent
2) Actually instead of $``\beth_\omega(\kappa) \le \mu"$ we just need a
conclusion of it: 
$$
(\forall \lambda \ge \mu)(\exists \theta)[\theta \in \text{{\rm Reg\/} }
 \and \kappa^\ast  \le  \theta  \leq  \mu  \and \text{{\rm cov\/}}
(\lambda,\theta^+,\theta^+,\theta) = \lambda],
\tag"{$(*)_1=(*)^1_{\mu,\kappa^\ast}$}"
$$
\medskip
\noindent
or even just a conclusion of that:
\bigskip
$$
\text{for every } \lambda \ge \mu \text{ for some } \theta < \mu,
\theta \ge \kappa^\ast \text{ we have: }\tag"{$(*)_2 = (*)^2_{\mu,\kappa^*}$}"
$$
\roster
\item "{{}}" $\otimes^\theta_\lambda = \otimes^{\theta,\kappa^\ast}_\lambda$:
  there is no family  ${\Cal P}$ of $> \lambda$ subsets of $\lambda$
  each of \newline
  cardinality $\theta$ with the intersection of any two having 
cardinality $<\kappa^\ast$.
\endroster
\endproclaim
\bigskip

\remark{\stag{3.6A} Remark}  1) The holding of $(\ast)_2$ is characterized in
\cite[\S6]{Sh:410}. \newline
2)  On earlier results concerning such problems and earlier 
history see Hajnal, Juhasz, Shelah \cite{HJSh:249}.
In particular, the following is quite a well known problem:
\bigskip
\roster
\item "{$\bigoplus$}"  Arhangelskii's problem: Can every topological 
space be divided into two pieces, such that no part contains
a closed homeomorphic copy of  ${}^\omega 2$  (or any topological 
space $Y$ such that every scattered set is countable, and the closure 
of a non-scattered set has cardinality continuum)?
\endroster
\bigskip
\noindent
3)  Note that the condition in $(\ast)_2$ holds if $\mu = 2^{\aleph_0} >
\aleph_\omega$,
$\kappa =\aleph_0$, $\kappa^\ast = \aleph_1$ and $\otimes^1_{\aleph_\omega}$
(from \scite{2.1}) (which holds e.g. if $V = V^P_0$, $P$ a c.c.c. 
forcing making the continuum  $> \beth^{V_0}_\omega$).  So in this case 
the answer to $\oplus$ is positive. \newline
\noindent
4)  Also if $\mu = 2^{\aleph_0} > \theta \ge \aleph_1$, and
$(\forall \lambda)[\lambda \ge 2^{\aleph_0} \Rightarrow
\otimes^{\theta,\aleph_1}_\lambda]$ then the answer to $\oplus$ in (2)
is yes; now on  $\otimes^{\theta,\aleph_1}_\lambda$ see 
\cite[\S6]{Sh:410}.\newline
\endremark
\bigskip

\demo{Proof}  We prove by induction on $\lambda \in [\mu,|X|]$ that:
\bigskip
\roster
\item "{$(*)_\lambda$}"  if $Z,Y$ are disjoint subsets of  $X$,
  $|Y| \le \lambda$, then there is a set \newline
  $Y^+$, $Y \subseteq Y^+ \subseteq X \backslash Z$,  $|Y^+| \le \lambda$
  and a function  $h:Y^+ \rightarrow \mu$ \newline
  such that: if  $A \in {\Cal P}^\ast$, $\kappa^\ast \le \theta < \mu$,
  $\otimes^\theta_\lambda$, $|A \cap Y^+| \ge \theta$ \newline
  and $|A \cap Z| < \mu$ then $h \restriction (A \cap Y^+)$ is onto $\mu$.
\endroster
\enddemo
\bigskip

\demo{Case 1}  $\lambda = \mu$, so $|Y| \le \mu$. \newline
\noindent
Without loss of generality $[B \subseteq Y \and |B| \le \kappa \and
|c\ell(B)| = \mu \Rightarrow c\ell(B) \backslash Z \subseteq Y]$.  Now
just note that ${\Cal P}_Y =: \{ c\ell(B) \cap Y:B \subseteq Y,|B| \le
\kappa,|c\ell(B) \cap Y| = \mu\}$ has cardinality $\le \mu = \mu^\kappa$,
and by the definition of ${\Cal P}^*$ (using the demand $|A \cap Z| <
\mu$ in $(*)_\lambda$), it suffices that $h$ satisfies: $[A \in {\Cal P}_Y 
\Rightarrow h \restriction Z$ is onto $\mu]$, which is easily accomplished.
\enddemo
\bigskip

\demo{Case 2}  $\lambda > \mu$. \newline

Let $\chi = \left( 2^\lambda \right)^+$, $\langle N_i:i \le \lambda \rangle$
an increasing continuous sequence of elementary submodels of
$({\Cal H}(\chi),\in,<^\ast_\chi)$, $\langle X,c\ell,Y,Z,\lambda \rangle \in  
N_0$, $\mu + 1 \subseteq N_0$, $\langle N_i:i \le j \rangle \in N_{j+1}$
(when $i < \lambda$) and $\| N_i \| = \mu + |i|$. \newline
\noindent
We define by induction on $i < \lambda$, a set  $Y^+_i$ and a function  
$h_i$ as follows: \newline
\smallskip
$(Y^+_i,h_i)$ is the $<^\ast_\chi$-first pair $(Y^\ast,h^\ast)$ such 
that:
\roster
\item "{(a)}"  $Y^\ast \subseteq X \setminus (Z \cup \dsize \bigcup
_{j < i}Y^+_j)$
\sn
\item "{(b)}"  $Y \cap N_i \backslash \dsize \bigcup_{j < i} Y^+_j 
\backslash Z \subseteq X \cap N_i \backslash \dsize \bigcup_{j < i} 
Y^+_j \backslash Z \subseteq Y^\ast$
\sn
\item "{(c)}"  $|Y^\ast| = \mu + |i|$
\sn
\item "{(d)}"  $h^\ast:Y^\ast \rightarrow \mu $
\sn
\item "{(e)}"  if  $A \in {\Cal P}^\ast $,  
$\otimes^\theta_{\mu +|i|}$, $\kappa^\ast \le \theta < \mu$,
  $|A \cap Y^*| \ge \theta$, \newline
  $|A \cap (Z \cup \dsize \bigcup_{j < i} Y^+_j) | < \mu$
  then  $h^* \restriction (A \cap Y^*)$ is onto $\mu$.
\endroster 
\medskip
Note:  $(Y^+_i,h_i)$  exists by the induction hypothesis applied to the 
cardinal  $\mu + |i|$  and the sets  $Z \cup \dsize \bigcup_{j < i}Y^+_j$,
$X \cap  N_i \backslash \dsize \bigcup_{j < i}Y^+_j$.  Also it is easy to
check that  $\langle (Y^+_j,h_j):j \le i\rangle \in N_{i+1}$ (as we always
choose ``the $<^\ast_\chi$-first", hence $Y^+_i \subseteq N_{i+1}$). \newline
\noindent
Let  $Y^+ = \dsize \bigcup_{i < \lambda} Y^+_i$,  $h = \dsize \bigcup_
{i < \lambda} h_i$.  Clearly $Y \subseteq \dbcu_{i < \lambda} N_i$ hence by 
requirement (b) clearly  $Y \subseteq Y^+$, (and
even $X \cap N_\lambda \backslash Z \subseteq Y^+)$; by requirements (c)
(and (a)) clearly  $|Y^+| = \lambda$, by requirement (a) clearly
$Y^+ \subseteq X \backslash Z$ and even $Y^+ = X \cap N_\lambda \backslash Z$.
\newline
\noindent
By requirements (a) + (d),  $h$  is a function from  $Y^+$ to  $\mu$.  Now 
suppose  $A \in {\Cal P}^\ast$,  $\otimes^\theta_\lambda$,
$\kappa^\ast \le \theta < \mu$,  $|A \cap Y^+| \ge \theta$, $|A \cap Z|
< \mu$;  we should prove  $``h \restriction (A \cap Y^+)$ is onto $\mu$".
So  $|A \cap N_\lambda| \ge \theta$.  Choose $(\delta^\ast,\theta^\ast)$
a pair such that: 
\roster
\widestnumber\item{(iii)}
\item "{(i)}"  $\delta^\ast \le \lambda$
\sn
\item "{(ii)}"  $\otimes^{\theta^*}_{\mu+|\delta^\ast|}$,  
$\kappa^\ast \le \theta^\ast < \mu$
\sn
\item "{(iii)}"  $|A \cap N_{\delta^\ast}| \ge \mu$ or $\delta^\ast =
\lambda$
\sn
\item "{(iv)}"  under (i) + (ii) + (iii), $\delta^\ast$ is minimal.
\endroster
\medskip
\noindent
This pair is well defined as  $(\lambda,\theta)$  satisfies requirement
(i) + (ii) + (iii).
\enddemo
\bigskip

\demo{Subcase 1}  $\delta^\ast$ is zero. 

So  $|Y^+_0 \cap  A| \ge \theta^\ast \ge \kappa^\ast$ hence by the 
choice of  $h_0$ we are done.
\enddemo
\bigskip

\demo{Subcase 2}  $\delta^\ast = i + 1$. 

So  $|A \cap  N_i| < \mu$,  hence  $|A \cap \dsize \bigcup_{j < i}Y^+_j| <
\mu$,  hence  $|A \cap (Z \cup \dsize \bigcup_{j < i}Y^+_j)| < \mu$.
Clearly $\otimes^{\theta^*}_{\mu +|i|}$ holds (as $\mu + |i| = 
\mu + |\delta^\ast|)$,  so if  $|A \cap Y^+_i| \ge \theta^\ast$ we are
done by the choice of  $h_i$; if not  $|A \cap (Z \cup \dsize \bigcup
_{j < i+1}Y^+_j)| < \mu$  and  $A \cap  Y^+_{i+1} \supseteq A \cap
N_{i+1} = A \cap N_{\delta ^\ast}$ has cardinality  $\ge \theta^\ast$
(and $\otimes^{\theta^\ast}_{|Y^+_{i+1}|}$ holds)
so we are done by the choice of  $h_{i+1}$.
\enddemo
\bigskip

\demo{Subcase 3}  $\delta^\ast$ limit  $< \lambda$. \newline
\noindent
So for some  $i < \delta^\ast$,  $|A \cap N_i| \ge \theta^\ast$
[why? as $\theta^\ast < \mu < \lambda$].  Now in
$N_{i+1}$ there is a maximal family  $Q \subseteq [X \cap N_i]^{\theta^\ast}$
satisfying  $[B_1 \ne B_2 \in Q \Rightarrow  |B_1 \cap B_2| < \kappa^\ast]$
hence  $|Q| \le \mu + |\delta^\ast|$ and without loss of generality $Q \in
N_{i+1}$,  hence $Q \subseteq  N_{\delta^\ast}$ so there is  $B \in  Q$,
$B \in N_{\delta^\ast}$,  $|B \cap  A| \ge \kappa^\ast$,  but $|B| =
\theta^\ast < \mu  = \mu^\kappa$ 
hence  $\left[ B^\prime \in  [B \cap A]^\kappa \Rightarrow  B \cap A \in  
N_{\delta^\ast} \right] $.  As  $A \in {\Cal P}^\ast$ there is $B^\prime
\in [B \cap A]^\kappa$ with $c\ell(B') \subseteq A$, $|c \ell(B')|
 = \mu$.  Clearly $c \ell(B') \in  N_{\delta^\ast}$ hence for some
$j \in (i,\delta^\ast)$,  $c \ell(B') \in N_j$ hence $c \ell(B') \subseteq 
X \cap N_j$.  So  $|A \cap N_j| \ge \mu$.  By assumption for some
$\theta^\prime \in [\kappa^\ast,\mu)$,  $\otimes^{\theta^\prime}_{\mu +|j|}$,
so $(j,\theta^\prime)$  contradicts the choice of $(\delta^\ast,\theta^\ast)$.
\enddemo
\bigskip

\demo{Subcase 4}  $\delta^\ast$ limit $=\lambda$. \newline
\noindent
As  $\lambda \in N_0$,  there is a maximal family  $Q \subseteq  
[\lambda]^{\theta^\ast}$ satisfying  $[B_1 \ne B_2 \in Q \Rightarrow  
|B_1 \cap B_2| < \kappa^\ast]$  which belongs to  $N_0$.  By the assumption 
$(*)_2$, we know  $|Q| \le \lambda$.  We define by induction on 
$j \le \lambda$,  a one-to-one function  $g_j$ from  $N_j \cap X \backslash
Z$  onto an initial segment of  $\lambda$ increasing continuous in $j$,
$g_j$ the $<^\ast_\chi$-first such function.  So clearly $g_j \in N_{j+1}$
and let  
$Q^\prime = \{ g_\lambda (B):B \in  Q \}$  (i.e. $\{ \{ g_\lambda
(x):x \in B \}: B \in Q \}$ note:  $g_\lambda$ is necessarily a one 
to one function from  $N_\lambda \cap X\backslash Z$  onto  $\lambda$).  So 
for some  $B \in  Q^\prime$,  $|B^\prime \cap  A| \ge \kappa^\ast$, so as
in subcase 3, for some  $B^\prime \in N_\lambda$,  $B^\prime \subseteq B \cap
A$,  $|B^\prime| =  \kappa,c \ell(B') \subseteq A,|c \ell(B')| = \mu $;  
so for some $i < \lambda,c \ell(B') \subseteq N_i$.
But $|A \cap Z| < \mu$ so $|A \cap Y^+_i| = \mu$  and by assumption
$(\ast)_2$, for some $\theta$, $\kappa^\ast \le \theta < \mu$  we have  
$\otimes^\theta_{\mu +|i|}$,  contradicting the choice of  
$(\delta^\ast,\theta^\ast)$ (i.e. minimality of  $\delta^\ast$).
\hfill$\square_{\scite{3.6}}$
\enddemo
\bigskip

\demo{\stag{3.7} Discussion}  1) So if 
we return to the topological problem (see
$\oplus$ of \scite{3.6A}(2)), by \scite{3.6} + \scite{3.6A}(4), 
if $2^{\aleph_0} > \theta \ge \aleph_1$ we can
try $\theta = \aleph_2$, $\kappa^\ast = \aleph_0$, $\kappa = \aleph_1$.
So a negative answer to $\oplus$ (i.e. the consistency of a negative
answer) is hard to come by: it implies that for some
$\lambda$, $\neg \otimes^{\theta,\aleph_1}_\lambda$, a statement which, when
$\theta > \aleph_1$ at present we do not know is consistent (but clearly it
requires large cardinals). \newline
2) If we want $\mu = 2^{\aleph_0} = \aleph_2$, $\theta = \aleph_1 =
\kappa^\ast$ we should consider a changed framework.  We have a family
${\frak I}$ of ideals on cardinals $\theta < \mu$ which are $\kappa$-based
(i.e. if $A \in I^+$, $I \in {\frak I}$ (similar to \cite{HJSh:249})
then $\exists B \in [A]^\kappa(B \in I^+$)) and in \scite{3.6} 
replace ${\Cal P}^\ast$ by

$$
\align
{\Cal P}^\ast = {\Cal P}^\ast_{\frak I} =: \biggl\{ A \subseteq X:&|A| = \mu
\text{ and for every pairwise distinct } \\
  &x_\alpha \in A \text{ for } \alpha < \theta \text{ we have } \\
  &\{ u \subseteq \theta:|c\ell \{x_\alpha:\alpha \in u \}| < \mu \} \\
  &\text{ is included in some } I \in {\frak I} \biggr\}.
\endalign
$$
and replace $(*)_2$ by
\roster
\item "{$(*)_3$}" For every $\lambda \ge \mu$ assume \newline
$F \subseteq
\{ (\theta,I,f):I \in {\frak I}, \theta = \text{ Dom}(I), f:\theta
\rightarrow \lambda \text{ is one to one}\}$
\endroster
and if $(\theta_\ell,I_\ell,f_\ell) \in F$ for $\ell = 1,2$ are distinct
then $\{ \alpha < \theta_2:f_2(\alpha) \in \text{ Rang } f_1\} \in I_2$.

Then $|F| \le \lambda$. \newline
Note that the present ${\Cal P}^\ast$ fits for dealing with $\oplus$ of 
\scite{3.6A}(2) and repeating the proof of \scite{3.6}. 
\enddemo
\bn
\ub{\stag{3.8} Discussion of Consistency of no}:  There are some restrictions
on such theorems.  Suppose
\mr
\item "{$(*)$}"  GCH and there is a stationary $S \subseteq \{\delta <
\aleph_{\omega +1}:\text{cf}(\delta) = \aleph_1\}$ and \nl
$\langle A_\delta:
\delta \in S \rangle$ such that: $A_\delta \subseteq \delta = \sup A_\delta$,
otp$(A_\delta) = \omega_1$ and \nl
$\delta_1 \ne \delta_2 \Rightarrow
|A_{\delta_1} \cap A_{\delta_2}| < \aleph_0$.
\ermn
(This statement is consistent by \cite[4.6,p.384]{HJSh:249} which continues
in \cite{Sh:108}.) \nl
Now on $\aleph_{\omega_1}$ we define a closure operation:

$$
\alpha \in c \ell(u) \Leftrightarrow (\exists \delta \in S)[\alpha \in
A_\delta \and (u \cap A_\delta) \ge \aleph_0].
$$
\mn
This certainly falls under the statement of \scite{3.6}(2) with $\kappa =
\kappa^* = \aleph_0,\mu = \aleph_1$ except the pcf assumptions $(*)_1$ and
$(*)_2$ fail.  However, this is not a case of our theorem.
\bn
\centerline {$* \qquad * \qquad *$}
\newpage

\head {Appendix: Existence of tiny models} \endhead  \resetall \bigskip

We deal now with a model theoretic problem, the existence of tiny models;
we continue Laskowski, Pillay, and Rothmaler \cite{LaPiRo}; our main result is
in \scite{ap.6}.
\demo{\stag{ap.1} Context} Assume $T$ is a complete first order theory.  
Let $|T|$ be the number of first order formulas $\varphi(\bar x)$,
$\bar x = \langle x_\ell:\ell < n \rangle$, $n < \omega$, up to 
equivalence modulo $T$. \newline
Assume $T$ is categorical in all cardinals $\chi > \lambda =: |T|$ and 
call a model $M$ of $T$ \underbar{tiny} if $\| M \| < \mu(=|T|)$.  It is
known that a $T$ with a tiny model satisfies exactly one of the following:
\medskip
\roster
\item "{(a)}" $T$ is totally transcendental, trivial (i.e. any regular type
is trivial)
\sn
\item "{(b)}"  $T$ is not totally transcendental.
\endroster
\enddemo
\bigskip

\subhead{\stag{ap.2} Question} \endsubhead  For which $\mu < \lambda$ are there
$T$, $|T| = \lambda$ (which is categorical in $\lambda^+$ and)
with a tiny model of cardinality $\mu$?
\bigskip

\demo{\stag{ap.3} Discussion}  By \cite{LaPiRo} we can deal 
with just the following two cases \newline
(see \cite{LaPiRo}, 0.3,p.386 and 387$^{1-21}$ and 1.7,p.390).
\enddemo
\bigskip

\subhead{Case A} \endsubhead $x=x$  is a minimal formula and its prime model
consists of individual constants.
\bigskip

\subhead{Case B}  \endsubhead $T$ is superstable not totally transcendental
and is unidimensional, the formula $x=x$ is weakly minimal, regular types 
are trivial and its prime model consists of individual constants.

They proved: $(\forall \kappa)[\kappa^{\aleph_0} \le \kappa^+ \Rightarrow
\text{ in case A, } \mu = \aleph_0]$, (see [LaPiRo,2.1,p.341]).  Actually more
is true by continuing their argument.
\bigskip

\proclaim{\stag{ap.4} Lemma} If $\lambda,\mu,T$ are as above, in Case A, 
\ub{then}:
\medskip
\roster
\item "{(i)}"  $\lambda < \beth_\omega$,
\sn
\item "{(ii)}"  we can find $\langle \lambda_n:n < \omega \rangle$ such that:
$\lambda_0 = \mu,\lambda_n \le \lambda_{n+1}$, \nl
$\lambda = \dsize \sum_{n < \omega} \lambda_n$ and
$(*)_{\mu,\lambda_n,\lambda_{n+1}}$ \nl
(hence in particular $(*)_{\mu,\mu,\mu^+}$), where
{\roster
\itemitem{ $(*)_{\mu,\sigma,\theta}$ }  there is a family of $\theta$ subsets
of $\sigma$ each of cardinality $\mu$, with the intersection of any two being
finite, or equivalently $\theta$ functions from $\mu$ to $\sigma$ such that
for any two such two distinct functions $f',f''$ we have 
$\{ i < \mu:f'(i) = f''(i)\}$ is finite.
\endroster}
\endroster
\endproclaim
\bigskip

\demo{Proof}  By \S1, (ii) $\Rightarrow$ (i), so let us
prove (ii).  Let $M$ be a tiny model of $T$, $\| M \| = \mu$.

For $n \ge 0$, let ${\frak B}_n$ be the family of definable (with parameters)
subsets of ${}^{n+1}M$.  Clearly
$|T| \le \dsize \sum_{n < \omega} |{\frak B}_n|$,
also $\mu = \| M \| \le |{\frak B}_n|$,
$|{\frak B}_n| \le |{\frak B}_{n+1}|$.  Also
$|{\frak B}_0| = \| M \|$ as $M$ is minimal which means $\lambda_0 = \mu$; let
$\lambda_n =: |{\frak B}_n|$, so $\lambda_n \le \lambda_n;\mu =
\dsize \sum_{n < \omega} \lambda_n$ and it is enough to prove 
$(*)_{\mu,\lambda_n,\lambda_{n+1}}$ when $\lambda_n < \lambda_{n+1}$.  
For each
$R \in {\frak B}_{n+1}$ we define a function $f_R$ from $M$ to ${\frak B}_n$,
$f_R(a) = \{ \bar b \in {}^n M:\bar b\char 94 <a> \in R \}$.  So
$\{f_R:R \in {\frak B}_{n+1} \}$ is a family of $\lambda_{n+1}$ functions
from $M$ to ${\frak B}_n$, hence it is enough to show:

$$
\text{define } R_1 \approx R_2 \Rightarrow \{ a \in M_0:f_{R_1}(a)
= f_{R_2}(a) \} \text{ is co-finite}
$$

\noindent
then
\roster
\item "{$(\alpha)$}"  $\approx$ is an equivalence relation on ${\frak B}
_{n+1}$
\sn
\item "{$(\beta)$}"  each $\approx$-equivalence class has cardinality
$\le \lambda_n$
\sn
\item "{$(\gamma)$}" if $\neg[R_1 \approx R_2]$, $R_1 \in {\frak B}_{n+1}$,
$R_2 \in {\frak B}_{n+1}$ then \nl
$\{ a \in M:f_{R_1}(a) = f_{R_2}(a) \}$ is finite.
\endroster
\medskip

Now clause $(\alpha)$ is straight, for clause $(\beta)$ use $x=x$ is
minimal and compute, for clause $(\gamma)$ remember
$x=x$ is a minimal formula.  Together, a set of representations $\varUpsilon$
for ${\frak B}_{n+1} / \approx$ will have cardinality $\lambda_{n+1}$ (as
$|{\frak B}_{n+1}| = \lambda_{n+1} > \lambda_n = |{\frak B}_n| \ge \mu$ 
by clauses $(\alpha),(\beta))$ and $\{f_R:R \in \varUpsilon \}$ is a set 
of functions as required.  \hfill$\square_{\scite{ap.4}\/}$
\enddemo
\bigskip

\proclaim{\stag{ap.5} Lemma}  Suppose $(*)_{\mu,\mu,\lambda}$, $\mu < \lambda$.
\ub{Then}
\medskip
\roster
\item "{(a)}"  there is a group $G$ of permutations of
$\mu$ such that $|G| = \lambda$ and \newline
$f \ne g \in G \Rightarrow \{ \alpha < \mu:f(\alpha) = g(\alpha) \}$ is
finite
\sn
\item "{(b)}"  there is a theory $T$ as in \scite{1.1}, 
$|T| = \lambda$, with a tiny model of cardinality $\mu$ of Case A.
\endroster
\endproclaim
\bigskip

\demo{Proof}  As (a) $\Rightarrow$ (b) is proved in 
\cite{LaPiRo}, p.392$^{23-31}$
we concentrate on (a).  Let $pr(-,-)$ be a pairing function on $\mu$ i.e.
$pr$ is one-to-one from $\mu \times \mu$ onto $\mu$.
So let $\{A_\zeta:\zeta < \lambda \} \subseteq
[\mu]^\mu$ be such that $\zeta \ne \xi \Rightarrow \aleph_0 > |A_\zeta \cap
A_\xi|$.  Clearly $\mu^{\aleph_0} \ge \lambda$ hence there is a list
$\langle \eta_\zeta:\zeta < \lambda \rangle$ of distinct members of
${}^\omega \mu$.  By renaming we can have the family $\{ A_{\zeta,n}:\zeta
< \lambda,n < \omega \}$, such that $(A_{\zeta,n} \in [\mu]^\mu,[(\zeta,n)
\ne (\xi,m) \Rightarrow |A_{\zeta,n} \cap A_{\xi,m}| < \aleph_0]$ and)
$\dsize \bigcup_{\zeta < \lambda} A_{\zeta,n} \cap
\dsize \bigcup_{\zeta < \lambda} A_{\zeta,m} = \emptyset$ for $n \ne m$, and
$\zeta \ne \xi \rightarrow (\exists n)(\forall m)[n \le m < \omega
\rightarrow A_{\zeta,n} \cap A_{\xi,n} = \emptyset]$.  
Let $g^0_{\zeta,n} \in {}^\mu \mu$ be $g^0_{\zeta,n}
(\alpha)$ = the $\alpha$th member of $A_{\zeta,n}$ and $g^1_{\zeta,n}(\alpha)
= pr(\alpha,g^0_{\zeta,n}(\alpha)$) so also $g^1_{\zeta,n}$ is a function
from $\mu$ to $\mu$.

We define the set $A = \mu \times ({}^{\omega >} \{-1,+1\} )$, clearly
$|A| = \mu;,$ let $x,y$ vary on $\{-1,+1\}$.  Now for $\zeta < \lambda$ we
define a permutation $f_\zeta$ of $A$, by defining $f^{+1}_\zeta \restriction
(\mu \times \{ \eta \}) = f_\zeta\restriction(\mu \times \{ \eta \})$,
$f^{-1}_\zeta\restriction(\mu \times \{ \eta \})$ for $\eta \in {}^n \{-1,+1
\}$ by induction on $n$ (so in the end, $f^{-1}_\zeta$ is the inverse of
$f_\zeta = f^{+1}_\zeta$). \newline
For $n=0, \eta = <>$ and let for $x \in \{ -1,+1 \}$, $f^x_\zeta(\alpha,<>) =
\langle g^1_{\zeta,0}(\alpha),<x> \rangle$.
\newline
For $n+1$, $\eta = \nu \char 94 <y> \, \in {}^{n+1} \{ -1,+1 \}$ we let
\medskip
\roster
\item "{$(\alpha)$}"  $f^x_\zeta(\alpha,\eta) = (\beta,\nu)$ \underbar{when}
\newline
$x = -y,f^y_\zeta(\beta,\nu) = (\alpha,\eta)$ (by the previous stage)
\sn
\item "{$(\beta)$}"  $f^x_\zeta(\alpha,\eta) = \langle g^1_{\zeta,n+1}
(\alpha),\eta\char 94 <x> \rangle$ \underbar{when} $(\alpha)$ does not apply.
\ermn
Easily $f_\zeta$ is a well-defined permutation of $A$. \nl
Now $\{ f_\zeta:\zeta < \lambda \}$ generates a group $G$ of permutations of
$A$.  We shall prove it generates $G$ freely, moreover:
\medskip
\roster
\item "{$\otimes$}"  \underbar{if} $n < \omega$, $t = \langle (\zeta(\ell),
x(\ell)):\ell \le n \rangle$ is such that $\zeta(\ell) < \lambda$, $x(\ell)
\in \{ -1,1 \}$, and for no $\ell < n$ do we have 
$\zeta(\ell) = \zeta(\ell + 1) \and x(\ell) = -x(\ell + 1)$ \nl
(i.e. $\dsize \prod_{\ell \le n} f^{x(\ell)}_{\zeta(\ell)}$ is a non-trivial
group term) \underbar{then} $A_t = \{ a \in A:(\dsize \prod_{\ell \le n}
f^{x(\ell)}_{\zeta(\ell)})(a) = a\}$ is finite.
\endroster
\medskip
\noindent
As $|A| = \mu$, this clearly suffices. \newline
As this property of $\dsize \prod_{\ell \le n} f^{x(\ell)}_{\zeta(\ell)}$
is preserved by conjugation \wilog \,
\mr
\item "{$(*)_0$}"  $\ell \le n \Rightarrow \zeta(\ell) \ne \zeta(\ell +1)
\vee x(\ell) \ne x(\ell +1)$ where $n+1$ is interpreted as zero.
\ermn
For any $a \in A_t$ let
\mr
\item "{$(*)_1$}"  $b^t_m[a] = (\dsize \prod^n_{\ell = m} f^{x(\ell)}
_{\zeta(\ell)})(a)$ for $m \le n + 1$ \newline
(so $b^t_{n+1}[a] = a = b^t_0[a]$ and for $m=0,\dotsc,n$ we have 
$b^t_m[a] = f^{x(m)}_{\zeta(m)}(b^t_{m+1}[a])$) 
\sn
\item "{$(*)_2$}"  $b^t_m[a] = \langle \beta^t_m[a],\eta^t_m[a] \rangle$
\ermn
Choose $m^* < \omega$ large enough such that:
\mr
\item "{$(*)_3$}"  if $m \ge m^*$ and $0 \le \ell_1 < \ell_2 \le n$ then
\newline
$A_{\zeta(\ell_1),m} \cap A_{\zeta(\ell_2),m} = \emptyset$.
\ermn
For $a \in A_t$ let $m = m[a] \le n+1$ be such that $lg(\eta^t_m[a])$ is
maximal and call the length $k=k[a]$.  
As $f_\zeta(\langle \alpha,\eta \rangle) = \langle \beta,\nu
\rangle$ implies $lg(\eta) \in \{ lg(\nu) -1,lg(\nu) + 1\}$, clearly
\mr
\item "{$(*)_4$}"  $lg(\eta^t_{m-1}[a]) = 
lg(\eta^t_{m+1}[a]) = lg(\eta^t_m[a]) - 1$ (where
$m-1,m+1$ means $\text{mod } n+1$). 
\ermn
Clearly
\mr
\widestnumber\item{$(*)_5(a)$}
\item "{$(*)_5(a)$}"  $b^t_m[a] = f^{x(m)}_{\zeta(m)}
(b^t_{m+1}[a])$
\sn
\item "{$(b)$}"  $b^t_{m-1}[a] = f^{x(m-1)}_{\zeta(m-1)}(b^t_m[a])$
hence (as $(f^{x(m-1)}_{\zeta(m-1)})^{-1} = f^{-x(m-1)}
_{\zeta(m-1)}$) we have
\sn
\item "{$(b)'$}"  $b^t_m[a] = f^{-x(m-1)}_{\zeta(m-1)}(b^t_{m-1}[a])$.
\ermn
Looking at the definition of $f^{-x(m-1)}_{\zeta(m-1)}(b^t_{m-1}[a])$, as
$m=m[a]$ by $(*)_4$ clause $(\beta)$ in the definition of $f$ applies so
\mr
\item "{$(*)_6(a)$}"  $f^{-x(m-1)}_{\zeta(m-1)}(b^t_{m-1}[a]) =
(g^1_{\zeta(m-1),k[a]}(\beta^t_{m-1}[a]),(\eta^t_{m-1}[a]) \char 94
\langle - x(m-1) \rangle)$. 
\ermn
Similarly looking at the definition $f^{x(m)}_{\zeta(m)}(b^t_{m+1}[a])$,
by $(*)_4$ clause $(\beta)$ applies so
\mr
\item "{$(*)_6(b)$}"  $f^{x(m)}_{\zeta(m)}(b^t_{m+1}[a]) =
(g^1_{\zeta(m),k[a]}(\beta^t_{m+1}[a]),(\eta^t_{m+1}[a]) \char 94
\langle x(m) \rangle)$.
\ermn
By $(*)_5(b)' + (*)_6(a)$ we have
\mr
\item "{$(*)_7(a)$}"  $b^t_m[a] =
(g^1_{\zeta(m-1),k[a]}(\beta^t_{m-1}[a]),(\eta^t_{m-1}[a]) \char 94
\langle - x(m-1) \rangle)$.
\ermn
By $(*)_5(a) + (*)_6(b)$ we have
\mr
\item "{$(*)_7(b)$}"  $b^t_m[a] =
(g^1_{\zeta(m),k[a]}(\beta^t_{m+1}[a]),(\eta^t_{m+1}[a]) \char 94
\langle x(m) \rangle))$.
\ermn
We can conclude by $(*)_7(a) + (*)_7(b)$ that
\mr
\item "{$(*)_8$}"  $x(m) = - x(m-1)$ hence $x(m) \ne x(m-1)$.
\ermn
So by $(*)_0$ applied to $m-1$ we get
\mr
\item "{$(*)_9$}"  $\zeta(m) \ne \zeta(m-1)$.
\ermn
Clearly by $(*)_7(a) + (*)_7(b)$
\mr
\item "{$(*)_{10}$}"  $g^1_{\zeta(m),k[a]}(\beta^t_{m+1}[a]) =
g^1_{\zeta(m-1),k[a]}(\beta^t_{m-1}[a])$.
\ermn
Now by the choice of the $g^1_\zeta$'s (and the pairing function) and
$(*)_{10}$
\mr
\item "{$(*)_{11}$}"  $\beta^t_{m+1}[a] = \beta^t_{m-1}[a]$ and \nl
$g^0_{\zeta(m),k[a]}(\beta^t_{m+1}[a]) = g^0_{\zeta(m-1),k[a]}
(\beta^t_{m-1})$.
\ermn
So by $(*)_{11}$ and the choice of the $g^0_\zeta$'s
\mr
\item "{$(*)_{12}$}"  $g^0_{\zeta(m),k[a]}(\beta^t_{m+1}[a]) =
g^0_{\zeta(m-1),k[a]}(\beta^t_{m-1}) \in A_{\zeta(m),k[a]} \cap
A_{\zeta(m-1),k[a]}$.
\ermn
If $k[a] > m^*$ we get a contradiction (by $(*)_3$), so remembering $m=m[a]$
necessarily $lg(\eta^t_{m[a]}[a]) \le m^* + 1$, hence by the choice of
$m[a]$ we have 
$\dsize \bigwedge_{\ell} lg(\eta^t_\ell[a]) \le m^*$. \newline
So $\{ \langle \eta^t_\ell[a]:\ell < n+1 \rangle:a \in A_t \}$ is finite,
hence it suffices to prove for each $\bar \eta \in {}^{n+1}\{-1,1\}$
the finiteness of

$$
A_{t,\bar \eta} = \{ a \in A_t:\langle \eta^t_\ell[a]:\ell < n + 1 \rangle
= \bar \eta \}
$$
\mn
for any given $\eta$. \nl
As for $a \in A_{t,\bar \eta}$ we have $\ell g(\eta^t_m[a]) \le m^*$ for
$\ell \le n+1$, it is enough to prove that for each
$\bar k = \langle k_\ell:\ell \le n \rangle$ the following set is finite:

$$
A_{t,\bar \eta,\bar k} =: \{a \in A_{t,\bar \eta}:\ell g(\eta^t_\ell[a]) =
k_\ell \text{ for } < n+1\}.
$$
\mn
Let $K(\bar k) = \{\ell \le n+1:k_\ell$ is $\ge k_{\ell-1},k_{\ell+1}\}$
(i.e. a local maximum).
\sn
For each $m \in K(\bar k)$, the arguments in 
$(*)_3 - (*)_{12}$ apply, so by $(*)_{11}$, if $a \in  A_{t,\bar \eta,\bar k}$
then the value $\ell g(\eta^t_m[a])$ is determined and
$g^0_{\zeta(m),k_m}(\beta^t_{m+1}[a]) \in A_{\zeta(m),k_m} \cap
A_{\zeta(m-1),k_m}$, but the latter is finite so we can fix
$g^0_{\zeta(m),k_m}(\beta^t_{m+1}[a]) = \gamma_m$ but
$g^1_{\zeta(m),k_m}(\beta^t_{m+1}[a])$ can be computed from $\gamma =
g^0_{\zeta(m),k_m}(\beta^t_{m+1}[a])$ and $(\zeta(m),k_m)$ i.e. as
$pr$(otp$(A_{\zeta(m),k_m} \cap \gamma),\gamma_m)$. \nl
But by $(*)_7(b)$ the latter is $\beta^t_m[a]$ and as $\eta^t_m[a] =
\eta_m$ the value of $b^t_m[a]$ is uniquely determined.  Similarly by 
induction we can compute the other $b^t_{m'}[a]$ for every $m'$, in 
particular $b^t_0[a] = a$, so we are done.
\hfill$\square_{\scite{ap.6}}$
\enddemo
\bigskip

\demo{\stag{ap.6} Conclusion}  For a cardinal $\mu$, the 
following are equivalent:
\mr
\item "{(a)}"  there is a $T$ as in \scite{ap.5}(b) (i.e. $T$ categorical in
$|T|^+,|T| > \mu)$, with a tiny model $M,\| M \| = \mu$ as in Case
A above)
\sn
\item "{(b)}"  $(*)_{\mu,\mu,\mu^+}$
\sn
\item "{(c)}"  there is a group $G$ of permutations of $\mu$, $|G| = \mu^+$
such that for $g \in G$, $\{ \alpha < \mu:g(\alpha) = \alpha \}$ is finite
or is $\mu$.
\endroster
\enddemo
\newpage
\shlhetal

\newpage
    
REFERENCES.  
\bibliographystyle{lit-plain}
\bibliography{lista,listb,listx,listf,liste}

\enddocument

\bye